# State-Estimators for Chemical Reaction Networks of Feinberg-Horn-Jackson Zero Deficiency Type*

Madalena Chaves†
Eduardo D. Sontag
Department of Mathematics
Rutgers University
New Brunswick, NJ 08903
E-mail: {madalena,sontag}@math.rutgers.edu

**Abstract**

This paper provides a necessary and sufficient condition for detectability, and an explicit construction of observers when this condition is satisfied, for chemical reaction networks of the Feinberg-Horn-Jackson zero deficiency type.

*Keywords:* observers, chemical reaction systems, detectability

## 1 Introduction

One of the most interesting questions in control theory is that of constructing observers. Observers compute estimates of the internal states of a dynamical system, using data provided by measurement probes or partial state information. For linear systems, Luenberger observers (also known as "deterministic Kalman filters" since they amount to Kalman filters designed without regard to the statistics of measurement noise) provide a general solution, but, for nonlinear systems, establishing generally applicable conditions for existence and convergence of observers is an open and active area of research.

The question of constructing observers has long been of interest in chemical engineering, and in particular for bioreactors; see for instance the papers [2, 5, 15] and the book [3]. Observers are of great potential interest, in particular, for biological experimentation and biomedical applications, where the online monitoring of proteins involved in signaling pathways (using for instance fluorescent labelings of molecules) will lead to a better understanding of cellular dynamical processes.

The main objective of this paper is to construct (when and if possible) state observers for the class of systems which describe chemical reaction networks of the type introduced by Feinberg, Horn, and Jackson in [6, 7, 8, 9, 11, 12]. As outputs, we take a subset of state variables or, more generally, monomials in state variables (which could, in practice, be associated to measured reaction rates). We provide here a complete solution to the observer problem for this class. The results given in [20] provide a convenient formalism as well as a set of technical tools, in the form of Lyapunov estimates, which are central to our results, and we will repeatedly refer to that paper for basic concepts and results.

*Work was supported in part by US Air Force Grant F49620-98-1-0242
†Supported by Fundação para a Ciência e a Tecnologia, Portugal, under the grant PraxisXXI/BD/11322/97

As a first step in this paper, we prove a necessary and sufficient theoretical condition for detectability. As a second step, we proceed to explicitly construct two full-state observers that are guaranteed to converge globally, under the hypothesis that the system is detectable.

In terms of performance, the second observer seems to have substantially slower convergence characteristics but, on the other hand, it is more "robust" to perturbations in measurements. This robustness can be quantified by an input-to-state stability (ISS) type of estimate, and we prove that the second observer is ISS with respect to a far wider class of external disturbances than the first (more efficient) observer.

We also provide simulations which test the behavior of our observers in the presence of observation noise and even of unknown inputs acting on the system; the observers turn out to be surprisingly robust to such effects, but we leave for future work the formulation of theoretical results which quantify this robustness. The advantadges of our observers over the standard constructions of Luenberger and the extended Kalman filter are also illustrated by simulations. Most of the examples worked out in this paper, as well as the simulations testing the robustness of our observers, all concern the kinetic proofreading model proposed by McKeithan in [14] for T-cell receptor signal transduction (which motivated [20]).

The organization of this paper is as follows:

- This section introduces the problem as well as some notations and definitions to be adopted in this work.

- Sections 2 and 3 introduce the observers, and the main results are stated and proved. Some examples that illustrate the result are presented, applying the construction to the case of the McKeithan network.

- Section 4 shows how to adapt the proof of the main result to a more general model.

- Section 5 provides simulations, testing in particular the effect of noise and of unknown inputs acting on states, or state drift. A comparison between the Luenberger observer, an extended Kalman filter and our main observer is also given here.

- An appendix collects various technical results.

## 1.1 The Problem

The main objective of this paper is to construct (when and if possible) observers for a certain type of systems that provide a mathematical model for a class of chemical reactions. Several stability and other control-theoretic results applicable to such systems were discussed and surveyed (but without outputs) in the paper [20] (see also [4]). The systems that we study all have the generic form

$$\begin{aligned} \dot{x} &= f(x) \\ y &= h(x), \end{aligned} \qquad (1)$$

with the requirements on $f$ and $h$ specified next.

The function $f : \mathbb{R}^n \to \mathbb{R}^n$ is of the mass-action kinetics form

$$\sum_{i=1}^{m} \sum_{j=1}^{m} a_{ij} x_1^{b_{1j}} x_2^{b_{2j}} \ldots x_n^{b_{nj}} (b_i - b_j), \qquad (2)$$

where $m \leq n$ and each $b_j$ is a column vector in $\mathbb{R}^n$ and has entries $b_{1j}, b_{2j}, \ldots, b_{nj}$, which are nonnegative integers. It is assumed that $B := [b_1, b_2, \ldots, b_n]$ has rank $m$ and that none of its rows vanishes.

The constants $a_{ij}$ are all nonnegative, and the matrix $A = [a_{ij}]$ is assumed to be irreducible. In the language of Feinberg et. al., irreducibility amounts to a restriction to "single linkage class" systems. This restriction can be removed, as explained in [20], provided that the space $\mathcal{D}$ introduced below is defined in a slightly different way to account for the number of connected components in the incidence graph of $A$. In order to simplify the presentation, the main result is stated for irreducible systems, and a sketch of how to treat the "multiple linkage classes" case is given in Section 4.

Although (2) is defined on all of $\mathbb{R}^n$, we will be interested only on those trajectories which evolve in the positive orthant $\mathbb{R}^n_{>0} = \{(x_1, \ldots, x_n) \in \mathbb{R}^n : x_i > 0 \text{ for all } i\}$. It is easy to verify (cf. [20] and below) that $\mathbb{R}^n_{>0}$ is a forward invariant set for (1) when $f$ has the form (2).

Define
$$\mathcal{D} := \text{span}\{b_i - b_j : i, j = 1, \ldots, m\}$$
and let the canonical basis of $\mathbb{R}^n$ be $\{e_i : i = 1, \ldots, n\}$.

The function $f$ does not depend explicitly on time $t$; however, throughout the text, several variations of the system will be considered, obtained by adding input terms to $f$. The inputs will be assumed to be measurable and bounded functions $u : [0, +\infty) \to \mathbb{R}^p$, and the resulting right hand side will be denoted by $f^*(x, u)$.

In order to decide what kind of output functions $h : \mathbb{R}^n \to \mathbb{R}^p$ are natural to consider, we should think of the quantities that may be measured when performing a chemical experiment. Some possibilities are, for example, concentrations of some of the substances, or certain reaction rates (through markers, fluorescence, or energy released). This leads us to consider outputs whose coordinates are monomials. This kind of output includes both the case when the concentrations of some of the substances are measured ($x_1$, $x_2$, etc.) and the case when some of the reaction rates are measured (proportional to a monomial such as $x_1 x_2^3$).

Thus, we consider in this paper output maps $h : \mathbb{R}^n \to \mathbb{R}^p$ (typically, $p \leq n$), of the form:
$$h(x) = \begin{pmatrix} x_1^{c_{11}} x_2^{c_{12}} \cdots x_n^{c_{1n}} \\ \vdots \\ x_1^{c_{p1}} x_2^{c_{p2}} \ldots x_n^{c_{pn}} \end{pmatrix}, \tag{3}$$
where
$$C = \begin{bmatrix} c_{11} & c_{12} & \cdots & c_{1n} \\ c_{21} & c_{22} & \cdots & c_{2n} \\ & & \ddots & \\ c_{p1} & c_{p2} & \cdots & c_{pn} \end{bmatrix}$$
is a matrix all whose entries are either 0 or real numbers $\geq 1$. (In view of the preceding discussion, the most natural choice would be to take the entries of $C$ to be nonnegative integers, but we allow more arbitrary exponents since the results do not require integers. The restriction $c_{ij} \geq 1$ is imposed in order to insure that $h(x)$ is locally Lipschitz, which is needed in order to guarantee uniqueness of solutions in the observer equations to be presented later. Although we are ultimately interested in behavior for positive $x_i$'s, the ouputs make sense on $\mathbb{R}^n$, provided that we interpret exponents $x_i^c$ as $|x_i|^c$ for negative $x_i$'s.)

Let us introduce the following vector functions:
$$\rho_n(x) = (\ln x_1, \ldots, \ln x_n)'$$
defined on $\mathbb{R}^n_{>0}$, and
$$\text{Exp}_n(v) = (e^{v_1}, \ldots, e^{v_n})'$$

defined on $\mathbb{R}^n$. (From now on, we will drop the subscript $n$ of $\rho_n$ and $\text{Exp}_n$, since its value is usually clear from the context.) Then, for $x \in \mathbb{R}^n_{>0}$,

$$\rho(h(x)) = C\rho(x)$$

and

$$h(x) = \text{Exp}(C\rho(x)),$$

as long as all state variables (concentrations, when dealing with chemical models) $x_i$ are positive.

**No Boundary Equilibra Assumption**

For the rest of this paper, we will make the following assumption: *the system (1) has no boundary equilibrium in any positive stoichiometric class.* That is, if $x = (x_1, \ldots, x_n)$ is any vector with nonnegative components $x_i$, and some component $x_i$ of $x$ vanishes, and if $x - \bar{x} \in \text{span}\{b_i - b_j, \ i,j = 1, \ldots, m\}$ for some $\bar{x} \in \mathbb{R}^n_{>0}$, then $f(x) \neq 0$. This assumption amounts to saying that no reaction consistent with positive concentrations can be in equilibrium if one of the participating substances is at zero concentration. It is an assumption that is often satisfied in chemical reaction models, and is in particular satisfied in the main example (kinetic proofreading) to be discussed. It is possible to weaken this boundary assumption and still obtain significant (though more restricted) results (using the techniques developed in [20]) but we prefer not to do so in order to streamline the presentation.

Under the above assumption, the following result is a simple consequence of the theory developed by Feinberg et. al., see [20]. A (positive) class is any intersection of $\mathbb{R}^n_{>0}$ with an affine manifold of the form $\xi + \mathcal{D}$, where $\xi \in \mathbb{R}^n_{>0}$ and $\mathcal{D} = \text{span}\{b_i - b_j : i,j = 1, \ldots, m\}$. We use the notation $|x|$ for Euclidean norms.

**Theorem 1** *For each positive class $\mathcal{C}$ there exists a (unique) state $\bar{x} = \bar{x}_\mathcal{C} \in \mathbb{R}^n_{>0}$ which is a globally asymptotically stable point relative to $\mathcal{C}$, i.e., for each $x_0 \in \mathcal{C}$, the solution of $\dot{x} = f(x)$, $x(0) = x_0$ is defined for all $t \geq 0$, and $x(t) \to \bar{x}$ as $t \to \infty$, and for all $\varepsilon > 0$ there exists $\delta > 0$ such that, if $|\bar{x} - x_0| < \delta$, then $|\bar{x} - x(t)| < \varepsilon$ for all $t > 0$.*

## 1.2 Detectability

We will adopt the following

**Definition 1.1** *The system (1) is detectable (on $\mathbb{R}^n_{>0}$) if, for every two trajectories $x(\cdot)$ and $z(\cdot)$ in $\mathbb{R}^n_{>0}$, defined for all $t \geq 0$,*

$$h(x(t)) \equiv h(z(t)) \ \Rightarrow \ |x(t) - z(t)| \to 0 \text{ as } t \to \infty.$$

**Remark 1.2** Although we use it in this paper, this is not the most natural definition of detectability, because it is not "well posed" enough. In principle, one would want the definition of detectability to also include the property: "$h(x(t)) \approx h(z(t))$ for all $t$ implies $|x(t) - z(t)|$ is asymptotically near zero as $t \to \infty$", which can be formulated as an "incremental output to state stability" (or more generally, "incremental input/output to state stability", if there are inputs) property. Such a more general concept, in the style of [13] and [23], can also be studied, and a future paper will deal with that matter.

**Definition 1.3** By a (full-state) *observer* for (1) we mean a system $\dot{z} = g(z, h(x))$ evolving on a state space $\mathcal{X}$ which is an open subset of $\mathbb{R}^n$ containing $\mathbb{R}^n_{>0}$, such that, for each $x(0), z(0)$ in $\mathbb{R}^n_{>0}$, the composite system has solutions defined for all $t > 0$, and $|z(t) - x(t)| \to 0$ as $t \to +\infty$.

This is a weak definition, on that only "attraction" and not stability is required; however, in our proofs we achieve a stability property as well, as will follow from "KL" estimates. Detectability is, obviously, necessary for the existence of an observer.

The next observation is valid, in general, for any system $\dot{x} = f(x)$, $y = h(x)$ for which every trajectory converges to an equilibrium. For such a system, we denote by $E$ the set of all equilibria, i.e., all $\bar{x}$ such that $f(\bar{x}) = 0$.

**Lemma 1.4** Detectability is equivalent to:
$$[h(\bar{x}) = h(\bar{z}) \ \& \ \bar{x}, \bar{z} \in E] \ \Rightarrow \ \bar{x} = \bar{z}. \tag{4}$$

*Proof.* [necessity] Suppose that the system is detectable, and pick $\bar{x} \neq \bar{z}$ distinct elements of $E$ so that $h(\bar{x}) = h(\bar{z})$. Then $x(t) \equiv \bar{x}$ and $z(t) \equiv \bar{z}$ are two trajectories with $h(x(t)) \equiv h(z(t))$ and distinct limits, a contradiction.

[sufficiency] Suppose that (4) holds and pick any two trajectories with $h(x(t)) \equiv h(z(t))$. Since $h$ is continuous, this implies $h(\bar{x}) = h(\bar{z})$ for the limits of $x$ and $z$, which exist and belong to $E$ by the assumption. Then (4) says $\bar{x} = \bar{z}$, as we wanted to prove. ∎

This lemma applies, in particular, to our system (1) (with $f$ as in (2)), when we view (1) as restricted to the invariant set $\mathbb{R}^n_{>0}$. In order to be consistent with the notations in [20], we denote the equilibrium set as "$E_+$", emphasizing that we are only considering equilibria in the positive orthant.

Since the function Exp is one to one, for positive vectors $x, z \in \mathbb{R}^n_{>0}$ we have
$$h(x) = h(z) \ \Leftrightarrow \ C\rho(x) = C\rho(z),$$
so that the condition $h(\bar{x}) = h(\bar{z})$ becomes just
$$\rho(\bar{x}) - \rho(\bar{z}) \in \ker C.$$

Recall also this fact from [20]: if $\bar{x} \in E_+$, then, for any $\bar{z} \in \mathbb{R}^n_{>0}$,
$$\rho(\bar{x}) - \rho(\bar{z}) \in \mathcal{D}^\perp \ \iff \ \bar{z} \in E_+. \tag{5}$$

**Theorem 2** The following statements are equivalent:
  (a) The system (1) with $f$ as in (2) and $h$ as in (3) is detectable;
  (b) $\forall x, z \in \mathbb{R}^n_{>0}$, if $\rho(x) - \rho(z) \in \ker C$ and $x, z \in E_+$, then $x = z$;
  (c) $\mathcal{D}^\perp \bigcap \ker C = \{0\}$;
  (d) $\mathcal{D} + \mathrm{im}\, C' = \mathbb{R}^n$.

*Proof.* [(a) $\Leftrightarrow$ (b)] That condition (4) is equivalent to (b) follows immediately from the discussion above.

[(b) $\Rightarrow$ (c)] Pick any $y \in \mathcal{D}^\perp \bigcap \ker C$; we need to show that $y = 0$. Let $\bar{x}$ be any point of $E_+$ and put $\tilde{y} = \rho(\bar{x}) - y$, $\tilde{y} \in \mathbb{R}^n$. Then find $z = \mathrm{Exp}(\tilde{y}) \in \mathbb{R}^n_{>0}$ so that $\tilde{y} = \rho(z)$. Thus $y = \rho(\bar{x}) - \rho(z)$ with $\bar{x} \in E_+$ and $z \in \mathbb{R}^n_{>0}$. By definition of $y$, $\rho(\bar{x}) - \rho(z)$ is contained both in $\mathcal{D}^\perp$ and in $\ker C$. Condition (5) now implies that $z \in E_+$. By assumption (b), we now conclude that $\bar{x} = z$, or equivalently, $y = 0$ as wanted.

[(c) $\Rightarrow$ (b)] Let $x, z \in \mathbb{R}^n_{>0}$ satisfy both $\rho(x) - \rho(z) \in \ker C$ and $x, z \in E_+$. Then, from (5), if follows that $\rho(x) - \rho(z) \in \mathcal{D}^\perp$. Therefore, $\rho(x) - \rho(z) \in \mathcal{D}^\perp \bigcap \ker C$. By assumption (c) $\rho(x) - \rho(z) = 0$, and therefore, since $\rho(\cdot)$ is a bijective function on $\mathbb{R}^n_{>0}$, we have $x = z$.

[(c) $\Leftrightarrow$ (d)] This equivalence follows by duality. ∎

**Remark 1.5** Note that $\dim \mathcal{D} = m - 1$ and that detectability implies $m - 1 + \operatorname{rank} C' = n$. In the case the matrix $C$ has full rank, then $\operatorname{rank} C' = p$ and detectability implies $m - 1 + p = n$.

## 2 Constructing Observers

It has been shown in [20] that the system (1) can be steered to any of its equilibrium points in $E_+$ by applying additive error feedback. Namely, for any specific element $\bar{x} \in E_+$, a set $K$ of indices such that $\mathcal{D} + \operatorname{span}\{e_k : k \in K\} = \mathbb{R}^n$ and any positive constants $\gamma_k, k \in K$, the following system satisfies $\lim_{t \to \infty} |x(t) - \bar{x}| = 0$ for any initial conditions $x(0) \in \mathbb{R}^n_{>0}$:

$$\dot{x} = f(x) + \sum_{k \in K} \gamma_k (\bar{x}_k - x_k) e_k. \tag{6}$$

So, a natural starting point in constructing an observer for $\dot{x} = f(x)$, would be to take a copy of it and add a linear feedback term similar to the above:

$$\dot{z} = f(z) + \sum_{k \in K} \gamma_k (x_k - z_k) e_k.$$

We know that eventually the trajectories of the system we want to observe will satisfy $x(t) \to \bar{x}$, for some $\bar{x} \in E_+$, so it seems reasonable to expect that the dynamics of the $z$-system will become very close to the dynamics of the modified $x$-system, (6), as time tends to infinity.

To construct such an observer, we need the "right output", that is, a function $h(x)$ that allows the construction of the suitable $\gamma_k(x_k - z_k)e_k$ terms. Thus, this approach will only work when the output map is linear and of a very special form. For instance, if there exists a matrix $L$ such that

$$Ly = LCx = \sum_{k \in K} x_k e_k$$

where $K$ is such that $\mathcal{D} + \operatorname{span}\{e_k : k \in K\} = \mathbb{R}^n$ (for detectability), then the observer suggested above can be constructed (using $Ly$ instead of $y$ as input to the observer).

But in the general case such a matrix $L$ does not exist and a different, but related, construction will be used. The main result is:

**Theorem 3** Consider the system (1) and assume that it is detectable. Then the following system, with state space $\mathcal{X} = \mathbb{R}^n$, is an observer for the system (1):

$$\dot{z} = f(z) + C'(h(x) - h(z)). \tag{7}$$

The proof is given in Section 2.2.

**Example 2.1** Consider the system with $n = 4$ and $m = 3$, determined by the vectors:

$$b_1 = (1,1,0,0)', \quad b_2 = (0,0,1,0)', \quad \text{and} \quad b_3 = (0,0,0,1)'.$$

Then

$$\mathcal{D} = \text{span}\{(1,1,-1,0), (1,1,0,-1)\},$$

and for positive constants $k, k_3, k_4$ and $\beta_3$, the system (1) becomes the "McKeithan network",

$$\begin{aligned}
\dot{x}_1 &= -kx_1x_2 + k_3x_3 + k_4x_4 \\
\dot{x}_2 &= -kx_1x_2 + k_3x_3 + k_4x_4 \\
\dot{x}_3 &= kx_1x_2 - (k_3 + \beta_3)x_3 \\
\dot{x}_4 &= \beta_3x_3 - k_4x_4.
\end{aligned} \qquad (8)$$

I. Suppose that the output is given by $h(x) = (x_1, x_4)'$. We can take the following observer:

$$\begin{aligned}
\dot{z}_1 &= -kz_1z_2 + k_3z_3 + z_4z_4 \;\; +(x_1 - z_1) \\
\dot{z}_2 &= -kz_1z_2 + k_3z_3 + z_4z_4 \\
\dot{z}_3 &= kz_1z_2 - (k_3 + \beta_3)z_3 \\
\dot{z}_4 &= \beta_3z_3 - k_4z_4 \qquad +(x_4 - z_4).
\end{aligned}$$

II. Suppose that the output is given by $h(x) = (x_1x_2^2, x_1x_4)'$. Then, we can construct the following observer:

$$\begin{aligned}
\dot{z}_1 &= -kz_1z_2 + k_3z_3 + z_4z_4 \;\; +(x_1x_2^2 - z_1z_2^2) + (x_1x_4 - z_1z_4) \\
\dot{z}_2 &= -kz_1z_2 + k_3z_3 + z_4z_4 \;\; +2(x_1x_2^2 - z_1z_2^2) \\
\dot{z}_3 &= kz_1z_2 - (k_3 + \beta_3)z_3 \\
\dot{z}_4 &= \beta_3z_3 - k_4z_4 \qquad +(x_1x_4 - z_1z_4).
\end{aligned}$$

The remainder of this section will be devoted to proving this theorem. The basic idea is to study the stability properties of system (1) when a certain input is added to the function $f$, specifically, the system with right-hand side:

$$f^*(z, u) := f(z) + C'(u - h(z)). \qquad (9)$$

We will show that, for the system thus obtained, an "input to state stability" condition holds. The observer for (1) is obtained by letting the input be $u(t) = h(x(t))$.

The analysis of this system with inputs is interesting in its own right, since it provides a means of studying the behavior of the model under bounded inputs.

## 2.1 An ISS Property

The definition of an input-to-state stable (ISS) system was introduced in [18]. Here, we adapt this notion to deal with constrained inputs and relative equilibria, as well as positive states (i.e., those states with all coordinates in the strictly positive half-line).

From now on, whenever we mention an *input* $u(\cdot)$, we will mean a measurable essentially bounded function $u : [0, +\infty) \to \mathbb{R}^p$, possibly restricted to take values in a set $\mathbb{U}$ of $\mathbb{R}^p$. For $u : [0, +\infty) \to \mathbb{R}^p$ and any fixed $\bar{u} \in \mathbb{R}^n$, denote

$$\|u - \bar{u}\| := \text{ess.sup}.\{|u(t) - \bar{u}| : t \geq 0\}.$$

**Definition 2.2** A system $\dot{z} = f^*(z, u)$, with input-value set $\mathbb{U}$, evolving on a state space $\mathcal{X}$ which is an open subset of $\mathbb{R}^n$ containing $\mathbb{R}^n_{>0}$, is $\mathbb{R}^n_{>0}$-*(forward) invariant* if, for each initial state $z(0) \in \mathbb{R}^n_{>0}$ and each $\mathbb{U}$-valued input $u(\cdot)$, the corresponding maximal solution of $\dot{z} = f^*(z, u)$ as a differential equation in $\mathcal{X}$, which is defined on an interval $J_{z(0),u} = [0, t_{\max})$, has values $z(t) \in \mathbb{R}^n_{>0}$ for all $t \in J_{z(0),u}$.

The system is $\mathbb{R}^n_{>0}$-*(forward) complete* if it is $\mathbb{R}^n$-(forward) invariant and, for each $z(0) \in \mathbb{R}^n_{>0}$ and $\mathbb{U}$-valued input $u(\cdot)$, $J_{z(0),u} = [0, +\infty)$.

For the following definitions, fix points $\bar{x} \in \mathbb{R}^n_{>0}$ and $\bar{u} \in \mathbb{U}$, where $\mathbb{U}$ is a subset of $\mathbb{R}^p$.

**Definition 2.3** A system $\dot{z} = f^*(z, u)$, evolving on a state space $\mathcal{X}$ which is an open subset of $\mathbb{R}^n$ containing $\mathbb{R}^n_{>0}$, is *input-to-state stable with input-value set* $\mathbb{U}$ (with respect to the point $\bar{x}$ and the input $\bar{u}$) if it is $\mathbb{R}^n_{>0}$-complete and there exist a function $\beta$ of class $\mathcal{KL}$ and a function $\varphi$ of class $\mathcal{K}_\infty$ such that, for each $\mathbb{U}$-valued input $u(\cdot)$, and each initial condition $z_0 \in \mathbb{R}^n_{>0}$, for all $t$ it holds that

$$|z(t) - \bar{x}| \leq \beta(|z_0 - \bar{x}|, t) + \varphi(\|u - \bar{u}\|). \tag{10}$$

To study the stability properties of the system $\dot{z} = f^*(z, u)$ a Lyapunov-type technique is used, and the following definition is needed.

**Definition 2.4** An *ISS-Lyapunov function with respect to the point $\bar{x}$ and input $\bar{u}$*, for the system $\dot{z} = f^*(z, u)$ with inputs in $\mathbb{U} \subseteq \mathbb{R}^p$, evolving on a state space $\mathcal{X}$ which is an open subset of $\mathbb{R}^n$ containing $\mathbb{R}^n_{>0}$, is a continuous function $V : \mathbb{R}^n_{\geq 0} \to \mathbb{R}$, whose restriction to $\mathbb{R}^n_{>0}$ is continuously differentiable, which satisfies:

(i) For $z \in \mathbb{R}^n_{\geq 0}$, $V(z) \geq 0$ and $V(z) = 0 \Leftrightarrow z = \bar{x}$.

(ii) The set $\{z \in \mathbb{R}^n_{\geq 0} : V(z) \leq L\}$ is compact, for each positive constant $L$.

(iii) There exist two functions $\alpha, \gamma \in \mathcal{K}_\infty$ such that

$$\nabla V(z) \, f^*(z, u) \leq -\alpha(|z - \bar{x}|) + \gamma(|u - \bar{u}|)$$

for all $u \in \mathbb{U}$ and $z \in \mathbb{R}^n_{>0}$.

**Remark 2.5** This definition differs slightly from other definitions of "ISS-Lyapunov" functions, such as given in [21]. The difference is in the fact that here the function $V$ is only required to be differentiable in the set $\mathbb{R}^n_{>0}$, and it is not required to satisfy a decrease condition except at positive vectors. Observe that $V$ is not proper when restricted to the positive orthant: it remains finite as the boundary of $\mathbb{R}^n_{>0}$ is approached.

**Remark 2.6** For a function $V$ defined as above, there always exist $\mathcal{K}_\infty$ functions, $\nu_1, \nu_2$, such that

$$\nu_1(|z - \bar{x}|) \leq V(z) \leq \nu_2(|z - \bar{x}|) \tag{11}$$

for all $z \in \mathbb{R}^n_{\geq 0}$ (see Appendix A.1 for a construction of these bounds).

Next, we introduce our candidate ISS-Lyapunov function. Fix an $\bar{x} \in E_+$, and define the following function: $V : \mathbb{R}^n_{\geq 0} \to \mathbb{R}$:

$$V(z) = \sum_{i=1}^n \bar{x}_i \, g\left(\frac{z_i}{\bar{x}_i}\right) = \sum_{i=1}^n \bar{x}_i \left[\frac{z_i}{\bar{x}_i} \ln \frac{z_i}{\bar{x}_i} + 1 - \frac{z_i}{\bar{x}_i}\right] \qquad (12)$$

where $g : \mathbb{R}_{>0} \to \mathbb{R}$, $g(r) = r \ln r + 1 - r$, with the convention that $g(0) = 1$. The function $V$ is continuously differentiable on $\mathbb{R}^n_{>0}$ and has the following properties:

(i) It is positive definite on $\mathbb{R}^n_{\geq 0}$ with respect to $\bar{x}$, i.e., $V(z) \geq 0$ and $V(z) = 0 \Leftrightarrow z = \bar{x}$:

$$g'(r) = \ln r$$

so each $g(r)$ strictly increases on $(1, +\infty)$ and strictly decreases on $(0, 1)$. But since $g(1) = 0$, it must be that $g(r) > 0$ for all $r \in \mathbb{R}_{>0}$. As a sum of nonnegative quantities, $V(z)$ is also nonnegative and, moreover,

$$V(z) = 0 \Leftrightarrow g\left(\frac{z_i}{\bar{x}_i}\right) = 0 \,\forall\, i \Leftrightarrow z_i = \bar{x}_i \,\forall\, i,$$

that is, $V(z) = 0 \Leftrightarrow z = \bar{x}$.

(ii) For each constant $L > 0$ the set $\{z \in \mathbb{R}^n_{\geq 0} : V(z) \leq L\}$ is a compact subset of $\mathbb{R}^n_{\geq 0}$: $V(z) \leq L$ implies $g(z_i/\bar{x}_i) \leq L$ for all $i$; continuity of $g$ on $[0, +\infty)$ (since $\lim_{w \to 0} g(w) = 1$) implies that $z_i$ stays in a compact interval of $\mathbb{R}_{\geq 0}$ and hence $z$ stays in a compact subset of $\mathbb{R}^n_{\geq 0}$.

Finally, in Appendix A.1 examples of possible functions $\nu_1, \nu_2$ are given.

In the next lemma, we show that the existence of an ISS-Lyapunov function (according to Definition 2.4) for a given system, implies that the trajectories of that system satisfy an ISS estimate.

**Lemma 2.7** Consider an $\mathbb{R}^n_{>0}$-invariant system $\dot{z} = f^*(z, u)$, with input-value set $\mathbb{U}$. Suppose that either:

- the system is $\mathbb{R}^n_{>0}$-complete, or
- the state-space $\mathcal{X}$ contains $\mathbb{R}^n_{>0}$.

Fix a state $\bar{x} \in \mathbb{R}^n_{>0}$ and input value $\bar{u} \in \mathbb{U}$. Suppose that that there is some ISS-Lyapunov funtion $V$ with respect to $\bar{x}$ and $\bar{u}$.

Then, the system is input-to-state stable with input-value set $\mathbb{U}$ (with respect to the same $\bar{x}$ and $\bar{u}$).

*Proof.* This proof is very similar to what is done in the case of the usual definition of an ISS system (see [21], for instance).

According to the definition, $V$ satisfies an estimate of the form

$$\nabla V(z) \, f^*(z, u) \leq -\alpha(|z - \bar{x}|) + \gamma(|u - \bar{u}|),$$

for each $z \in \mathbb{R}^n_{>0}$ and each $u$ in the input value set $\mathbb{U}$, where $\alpha, \gamma \in \mathcal{K}_\infty$. On the other hand, we know from Remark 2.6 that there exist two class $\mathcal{K}_\infty$ functions, $\nu_1, \nu_2$, such that $\nu_1(|z - \bar{x}|) \leq V(z) \leq \nu_2(|z - \bar{x}|)$. Introducing a new $\mathcal{K}_\infty$ function $\chi$ given by $\chi(r) = \alpha^{-1}(2\gamma(r))$, observe that:

$$|z - \bar{x}| > \chi(|u - \bar{u}|) \Leftrightarrow \alpha(|z - \bar{x}|) > 2\gamma(|u - \bar{u}|) \Rightarrow \nabla V(z) \, f^*(z, u) < -\frac{1}{2}\alpha(|z - \bar{x}|). \qquad (13)$$

Now pick any initial condition $z_0 \in \mathbb{R}^n_{\geq 0}$, and an input $u : [0,+\infty) \to \mathbb{U}$, and consider the corresponding maximal solution $z(t)$ defined on $J = [0, t_{\max})$. For $s = \nu_2(\chi(\|u - \bar{u}\|))$, define the following sublevel set of $V$:

$$S = \{\xi \in \mathbb{R}^n_{\geq 0} : V(\xi) \leq s\}.$$

*Claim.* Suppose there exists an instant $\sigma \in J$ such that $z(\sigma) \in S$. Then $z(t) \in S$ for all $\sigma \leq t < t_{\max}$.

To see this, argue by contradiction: suppose there exists a $t > \sigma$ (but $t < t_{\max}$) and an $\varepsilon > 0$ such that $V(z(t)) > s + \varepsilon$. Let

$$\tau = \inf\{t \geq \sigma : V(z(t)) \geq s + \varepsilon\}.$$

Then $z(\tau) \notin S$ which implies $V(z(\tau)) > \nu_2(\chi(\|u - \bar{u}\|))$, and therefore, since $V(z) \leq \nu_2(|z - \bar{x}|)$ and $\nu_2$ is strictly increasing,

$$\nu_2(|z(\tau) - \bar{x}|) > \nu_2(\chi(\|u - \bar{u}\|)) \iff |z(\tau) - \bar{x}| > \chi(\|u - \bar{u}\|).$$

By (13), $\frac{d}{dt}V(z(t))|_\tau < 0$, implying that $V(z(t_*)) \geq V(z(\tau))$ for some $t_* \in (\sigma, \tau)$, and thus contradicting minimality of $\tau$. So the claim holds.

Now, let

$$T = \inf\{\sigma : z(\sigma) \in S\}$$

(with $T = +\infty$ if the trajectory never enters $S$). We have two cases to consider, for each $t \geq 0$:

For $t \geq T$: $V(z(t)) \leq \nu_2(\chi(\|u - \bar{u}\|))$ implies

$$\nu_1(|z(t) - \bar{x}|) \leq \nu_2(\chi(\|u - \bar{u}\|)) \implies |z(t) - \bar{x}| \leq \nu_1^{-1} \circ \nu_2 \circ \chi(\|u - \bar{u}\|),$$

so we let $\varphi = \nu_1^{-1} \circ \nu_2 \circ \chi$.

For $t < T$: $V(z(t)) > \nu_2(\chi(\|u - \bar{u}\|))$, which implies $|z(t) - \bar{x}| \geq \chi(\|u - \bar{u}\|)$ and hence

$$\frac{d}{dt}V(z(t)) \leq -\frac{1}{2}\alpha(|z(t) - \bar{x}|) \leq -\frac{1}{2}\alpha[\nu_2(V(z(t)))].$$

By a standard comparison principle, there exists a function $\tilde{\beta} \in \mathcal{KL}$ (which depends only on $\alpha$ and $\nu_2$) such that $V(z(t)) \leq \tilde{\beta}(V(z_0), t)$ for all $t < T$. Then

$$|z(t) - \bar{x}| \leq \nu_1^{-1}(\tilde{\beta}(V(z_0), t)) \leq \nu_1^{-1}(\tilde{\beta}(\nu_2(|z_0 - \bar{x}|), t)) := \beta(|z_0 - \bar{x}|, t).$$

Thus, for all $t \in J$,

$$|z(t) - \bar{x}| \leq \max\{\beta(|z_0 - \bar{x}|, t), \varphi(\|u - \bar{u}\|)\} \leq \beta(|z_0 - \bar{x}|, t) + \varphi(\|u - \bar{u}\|)$$

where $\beta \in \mathcal{KL}$ and $\varphi \in \mathcal{K}_\infty$.

Under the assumption that the system is $\mathbb{R}^n_{\geq 0}$-complete, we know that $J = [0, \infty)$, and thus this inequality holds for all $t \geq 0$.

If, instead, we assume that $\mathcal{X}$ contains $\mathbb{R}^n_{\geq 0}$, we argue as follows. The above ISS inequality implies that $\{z(t), t \in J\}$ is bounded (by $\beta(|z_0 - \bar{x}|, 0) + \varphi(\|u - \bar{u}\|)$) and is a subset of $\mathbb{R}^n_{\geq 0}$. Thus, the maximal solution remains in a compact set of $\mathcal{X}$ for all $t$ in the maximal interval $J$, which implies that $J$ has to be all of $[0, +\infty)$. ∎

## 2.2 Proof of Theorem 3

**Lemma 2.8** Let $f^*(z,u) = f(z) + C'(u - h(z))$ with $C$ such that $\mathcal{D} + \operatorname{im} C' = \mathbb{R}^n$. Let $\bar{x}$ denote any point in $E_+$, and $\theta$ be any real number with $0 < \theta < 1$. Define the following subset of $\mathbb{R}^p$:

$$\mathbb{U}_\theta = \{u \in \mathbb{R}^p : |u_k - h_k(\bar{x})| \leq \frac{\theta}{2} h_k(\bar{x}),\ k = 1, \ldots, p\}.$$

Then, for the function $V$ defined in (12), there exist functions $\alpha = \alpha_{\bar{x}}, \gamma = \gamma_{\bar{x}}$ of class $\mathcal{K}_\infty$ such that

$$\nabla V(z)\ f^*(z, u) \leq -\alpha(|z - \bar{x}|) + \gamma(|u - h(\bar{x})|),$$

for all $z \in \mathbb{R}^n_{>0}$ and all $u \in \mathbb{U}_\theta$.

*Proof.* Pick any $\bar{x} \in E_+$ and any $0 < \theta < 1$. We have

$$\nabla V(z)\ f^*(z, u) = \langle \rho(z) - \rho(\bar{x}), f(z) \rangle + \langle \rho(z) - \rho(\bar{x}), C'(u - h(z)) \rangle.$$

Using the notation $\varrho = \rho(z) - \rho(\bar{x})$, notice that the second term on the right can be rewritten as:

$$\langle C\varrho, u - h(z) \rangle = \langle C\varrho, u - h(\bar{x}) \rangle - \langle C\varrho, h(z) - h(\bar{x}) \rangle.$$

Introducing the notation $\sigma = C\varrho$, $\mu = h(z) - h(\bar{x})$ and $v = u - h(\bar{x})$, the expression for $\nabla V(z)\ f^*(z, u)$ becomes

$$\begin{aligned} \nabla V(z)\ f^*(z, u) &= \langle \varrho, f(z) \rangle - \langle \sigma, \mu \rangle + \langle \sigma, v \rangle \\ &= P(z, \bar{x}) + R(z, u, \bar{x}) \end{aligned}$$

where

$$P(z, \bar{x}) = \langle \varrho, f(z) \rangle - (1 - \theta)\langle \sigma, \mu \rangle$$

and

$$R(z, u, \bar{x}) = -\theta\langle \sigma, \mu \rangle + \langle \sigma, v \rangle = \sum \sigma_i(-\theta\mu_i + v_i)$$

We now bound each of these terms.

*Step 1.* We show first that $R(z, u, \bar{x}) \leq c_3|v|^2$, for some positive constant $c_3$.

Notice that $\mu_i \sigma_i = (h_i(z) - h_i(\bar{x}))(\ln h_i(z) - \ln h_i(\bar{x})) \geq 0$ for all pairs $h_i(z), h_i(\bar{x})$.

(i) if $\theta|\mu_i| \geq |v_i|$, then immediately

$$\sigma_i(-\theta\mu_i + v_i) \leq 0.$$

(ii) if $\theta|\mu_i| < |v_i|$, then

$$\sigma_i(-\theta\mu_i + v_i) \leq 2|\sigma_i||v_i| \leq \frac{2c_L}{\theta}|v_i|^2$$

where the last inequality follows from the bounds on $u$:

$$|v_i| = |u_i - h_i(\bar{x})| \leq \frac{\theta}{2} h_i(\bar{x}) \quad \Rightarrow \quad |h_i(z) - h_i(\bar{x})| = \mu_i \leq \frac{1}{\theta}|v_i| \leq \frac{1}{2} h_i(\bar{x}),$$

so that $h_i(z) \geq h_i(\bar{x})/2$, and with a Lipschitz constant, $c_L$, of the logarithmic function on $[A, +\infty)$, when $A = \min\{h_i(\bar{x})/2 : i = 1, \ldots, p\}$ (for instance, $c_L = 1/A$):

$$|\sigma_i| = |\ln h_i(z) - \ln h_i(\bar{x})| \leq c_L|h_i(z) - h_i(\bar{x})| = c_L|\mu_i|.$$

In either case,
$$R(z, u, \bar{x}) \leq \frac{2c_L}{\theta} \sum |v_i|^2$$
as wanted.

*Step 2.* Show that $P(z, \bar{x}) \leq -\alpha(|z - \bar{x}|)$, where $\alpha \in \mathcal{K}_\infty$.

Recall a result from [20], namely:
$$\begin{aligned}\langle \rho(z) - \rho(\bar{x}), f(z) \rangle &\leq -c(\bar{x}) \sum_{i=1}^{m} \sum_{j=1}^{m} \langle b_i - b_j, \rho(z) - \rho(\bar{x}) \rangle^2 \\ &\leq -c(\bar{x})|D_0\left(\rho(z) - \rho(\bar{x})\right)|^2\end{aligned} \quad (14)$$

where $D_0$ is an $(m-1) \times n$ matrix whose rows (transposed) constitute a basis of $\mathcal{D}$, and $c(\cdot)$ is a positive function. Using this estimate yields:
$$P(z, \bar{x}) \leq -c_1 |D_0 \varrho|^2 - (1 - \theta) \langle \sigma, \mu \rangle$$
where $c_1 = c(\bar{x})$. By assumption, $\mathcal{D} + \operatorname{im} C' = \mathbb{R}^n$, so we know that
$$\operatorname{span}\{\operatorname{cols}(D_0')\} + \operatorname{span}\{\operatorname{cols}(C')\} = \mathbb{R}^n,$$
so the matrix $D = [D_0' \ C']'$ has rank $n$.

Next, note that combining Lemmas A.2 and A.3, it is enough to find two $\mathcal{K}_\infty$ functions $\alpha_1, \alpha_2$ such that,
$$P(z, \bar{x}) \leq -\alpha_1(|D_0\varrho|) - \alpha_2(|C\varrho|).$$

We may simply take $\alpha_1(w) = c_1 w^2$ and will next produce $\alpha_2$.

Define, for each fixed $z$, the following finite, disjoint sets of integers
$$I_+ = I_+(z) = \{i : h_i(z) - h_i(\bar{x}) \geq 0\}$$
and
$$I_- = I_-(z) = \{i : h_i(z) - h_i(\bar{x}) < 0\}.$$

Write for each $i$,
$$\frac{h_i(z) - h_i(\bar{x})}{h_i(\bar{x})} = \ln h_i(z) - \ln h_i(\bar{x}) + g(h_i(z), h_i(\bar{x})).$$

It is not difficult to check that
$$g(h_i(z), h_i(\bar{x})) = \frac{h_i(z)}{h_i(\bar{x})} - \ln \frac{h_i(z)}{h_i(\bar{x})} - 1$$
has a minimum value of 0, achieved whenever $h_i(z) = h_i(\bar{x})$, and otherwise is strictly positive.

Let $C_i$ denote the $i$-th row of the matrix $C$ and let
$$\frac{\mu_i}{h_i(\bar{x})} = \xi_i + \zeta_i \quad \text{and} \quad \sigma_i = C_i\varrho = \xi_i + \eta_i$$
with
$$\xi_i = \begin{cases} \sigma_i = \ln h_i(z) - \ln h_i(\bar{x}) = C_i\varrho, & \text{if } i \in I_+ \\ \mu_i/h_i(\bar{x}) = (h_i(z) - h_i(\bar{x}))/h_i(\bar{x}), & \text{if } i \in I_- \end{cases}$$

$$\zeta_i = \begin{cases} g(h_i(z), h_i(\bar{x})), & \text{if } i \in I_+ \\ 0, & \text{if } i \in I_- \end{cases}$$

$$\eta_i = \begin{cases} 0, & \text{if } i \in I_+ \\ -g(h_i(z), h_i(\bar{x})), & \text{if } i \in I_-. \end{cases}$$

By construction, $(\xi_i + \zeta_i)(\xi_i + \eta_i) \geq \xi_i^2$, since all cross terms are nonnegative:

$$\xi_i \zeta_i = \begin{cases} g(h_i(z), h_i(\bar{x}))(\ln h_i(z) - \ln h_i(\bar{x})) \geq 0, & \text{if } i \in I_+ \\ 0, & \text{if } i \in I_- \end{cases}$$

$$\xi_i \eta_i = \begin{cases} 0, & \text{if } i \in I_+ \\ -g(h_i(z), h_i(\bar{x}))(h_i(z) - h_i(\bar{x}))/h_i(\bar{x}) \geq 0, & \text{if } i \in I_-, \end{cases}$$

and $\zeta_i \eta_i = 0$ for all $i$. Therefore, for each $i \in I_+$

$$\begin{aligned} \sigma_i \mu_i &= h_i(\bar{x})(\xi_i + \zeta_i)(\xi_i + \eta_i) \\ &\geq h_i(\bar{x}) \xi_i^2 \\ &= h_i(\bar{x})(C_i \varrho)^2, \end{aligned}$$

and for each $i \in I_-$

$$\begin{aligned} \sigma_i \mu_i &= h_i(\bar{x}) C_i \varrho \left( \frac{h_i(z) - h_i(\bar{x})}{h_i(\bar{x})} \right) \\ &= h_i(\bar{x}) |C_i \varrho| \left| \frac{h_i(z) - h_i(\bar{x})}{h_i(\bar{x})} \right| \\ &= h_i(\bar{x}) |C_i \varrho| \left| \frac{h_i(z)}{h_i(\bar{x})} - 1 \right| \\ &= h_i(\bar{x}) |C_i \varrho| \left| e^{\ln \left( \frac{h_i(z)}{h_i(\bar{x})} \right)} - 1 \right| \\ &= h_i(\bar{x}) |C_i \varrho| \left| e^{[\ln h_i(z) - \ln h_i(\bar{x})]} - 1 \right| \\ &= h_i(\bar{x}) |C_i \varrho| \left| e^{-|C_i \varrho|} - 1 \right| \end{aligned}$$

since $C_i \varrho = \ln h_i(z) - \ln h_i(\bar{x}) < 0$ for $i \in I_-$.

In summary, with $c_2 = \min\{h_1(\bar{x}), \ldots, h_p(\bar{x})\}$,

$$-(1-\theta)\langle \sigma, \mu \rangle \leq -c_2 \left( \sum_{i \in I_+} |C_i \varrho|^2 + \sum_{i \in I_-} |C_i \varrho| \left(1 - e^{-|C_i \varrho|}\right) \right),$$

and since both functions $a(r) = c_2 r^2$ and $b(r) = c_2 r(1 - e^{-r})$ are of class $\mathcal{K}_\infty$, one may apply Lemma A.3 repeatedly to obtain $\alpha_2 \in \mathcal{K}_\infty$ such that $-(1-\theta)\langle \sigma, \mu \rangle \leq -\alpha_2(|C\varrho|)$.

So, taking $\alpha \leq \alpha_1 + \alpha_2$ according to Lemma A.3, and $\gamma(r) = c_3 r^2$, gives the desired result. ∎

**Proposition 2.9** Suppose that the system defined by (2) and (3) is detectable. Consider the system with inputs

$$\dot{z} = f^*(z, u) := f(z) + C'(u - h(z)) \tag{15}$$

with state-space $\mathcal{X} = \mathbb{R}^n$. Then, the system is $\mathbb{R}^n_{>0}$-invariant with input-value set $\mathbb{R}^p_{\geq 0}$. Furthermore, let $\theta$ be any real number with $0 < \theta < 1$, and pick any fixed state $\bar{x} \in E_+$. Let $\mathbb{U}_\theta$ be the subset of $\mathbb{R}^p$ defined in Lemma 2.8. Then, the system (15) is ISS with input-value set $\mathbb{U}_\theta$ (with respect to the point $\bar{x}$ and the input $\bar{u} = h(\bar{x})$).

*Proof.* The proof of the first statement, namely that the system is $\mathbb{R}^n_{>0}$-invariant with input-value set $\mathbb{R}^p_{\geq 0}$, is fairly routine, and it proceeds as follows.

Given an initial condition $z(0) \in \mathbb{R}^n_{>0}$, and an $\mathbb{U}_\theta$-valued input, let $z(t)$ be the maximal solution of (15), defined on a (maximal) interval $J$. Let $\mathcal{I} = [0, +\infty)$.

Assume that one of the coordinates becomes $\leq 0$ at some instant and define

$$t_0 = \inf\{t \in J : z_k(t) = 0 \text{ for some } 1 \leq k \leq n\}.$$

Pick one coordinate $k$ such that $z_k(t_0) = 0$. We reorder variables, singling out this coordinate, and look at the time-dependent differential equation that results by fixing the remaining $n - 1$ variables. It is useful for that purpose to introduce the following notation:

$$(\check{z}(t), x) = (z_1(t), \ldots, z_{k-1}(t), x, z_{k+1}(t), \ldots, z_n(t)).$$

In addition, we wish to see the obtained scalar equation as well-defined for all $t$, not just $t \leq t_0$. So we construct a new function $F : \mathcal{I} \times \mathbb{R} \to \mathbb{R}$ as follows:

$$F(t, x) = \begin{cases} f^*_k(\check{z}(t), x; u(t)), & t \in [0, t_0) \\ f^*_k(\check{z}(t_0), x; u_0), & t \in [t_0, +\infty) \end{cases}$$

where $u_0$ is any fixed element of $\mathbb{U}_\theta$. Then, for each fixed $t$, $F(t, x)$ is locally Lipschitz in $x$ and the Lipschitz constants, $\alpha(t)$, are uniformly bounded (and hence locally integrable as a function of time). In addition, for each fixed $x$, $F(t, x)$ is measurable as a function of time. Thus the standard existence and uniqueness conditions apply.

*Claim.* $F(t, 0) \geq 0$ for almost all $t \in \mathcal{I}$.

To prove this, write

$$\begin{aligned} f^*_k(\check{z}, x; u) &= \sum_{i=1}^m \sum_{j \in A_0} a_{ij} z_1^{b_{1j}} \ldots z_{k-1}^{b_{(k-1)j}} z_{k+1}^{b_{(k+1)j}} \ldots z_n^{b_{nj}} b_{ki} \\ &+ \sum_{i=1}^m \sum_{j \in A_+} a_{ij} z_1^{b_{1j}} \ldots z_{k-1}^{b_{(k-1)j}} x^{b_{kj}} z_{k+1}^{b_{(k+1)j}} \ldots z_n^{b_{nj}} (b_{ki} - b_{kj}) \\ &+ \sum_{j=1}^p c_{jk}[u_j - h_j(\check{z}, x)], \end{aligned} \tag{16}$$

where $A_0 = \{j : b_{kj} = 0\}$ and $A_+ = \{j : b_{kj} > 0\}$.

For $x = 0$ and $t \in \mathcal{I}$:

(a) the third term is nonnegative since we are assuming that $c_{ji} \geq 0$ and $u_j \geq 0$ for all $i, j$, and because $h_j(\check{z}, x) = z_1^{c_{j1}} \ldots x^{c_{jk}} \ldots z_n^{c_{jn}}$, so either $c_{jk} = 0$, or $c_{jk} > 0$ and $h_j(\check{z}, 0) = 0$.

(b) the second term is zero since $x = 0$;

(c) the first term is nonnegative since, by definition of $t_0$, we are evaluating at $z_i = z_i(t) \geq 0$, for all $i$ and $t \leq t_0$, and $z_i = z_i(t_0)$ for $t > t_0$.

This proves the claim.

Moreover, notice that, for all $t \leq t_0$, the scalar variable $z_k(t)$ satisfies the initial value problem

$$\dot{x} = F(t, x)$$
$$x(0) = z_k(0),$$

where $F(t, 0) \geq 0$ for all $t \geq 0$. Solutions of this initial value problem exist on an open interval $\tilde{J}$, and this interval contains $[0, t_0]$ because $z_k(t)$ solves the equation in that interval. Then, by Lemma A.4, $x(t) > 0$ on $\tilde{J}$ and $z_k(t) = x(t) > 0$ for all $t < t_0$; since both $x(t)$ and $z_k(t)$ are continuous functions, we also have that $z_k(t_0) = x(t_0)$, contradicting the fact that $z_k(t_0) = 0$.

This concludes the proof of $\mathbb{R}^n_{>0}$-invariance with input-value set $\mathbb{R}^p_{\geq 0}$. As this implies that (15) is also $\mathbb{R}^n_{>0}$-invariant with (the smaller) input-value set $\mathbb{U}_\theta$, to prove that (15) is ISS with input-value set $\mathbb{U}_\theta$, it is enough, by Lemma 2.7, to show that this system admits the function $V$ defined in (12) as an ISS-Lyapunov function with respect to the state $\bar{x}$ and the input $h(\bar{x})$.

So we must check that $V$ satisfies Definition 2.4. Properties (i) and (ii) have already been shown in the discussion following formula (12). Property (iii) follows from Lemma 2.8. ∎

**Lemma 2.10** Let $f^*(z, u) = f(z) + C'(u - h(z))$ and let $V$ denote the function defined in (12). For each fixed $\bar{x} \in E_+$, and each constant $u_{\max} \geq 0$, there exists a constant $c_{\bar{x}, u_{\max}}$ such that

$$\nabla V(z)\, f^*(z, u) \leq c_{\bar{x}, u_{\max}}, \quad \forall\, z \in \mathbb{R}^n_{>0}, \quad \forall\, u \in [0, u_{\max}]^p.$$

*Proof.* Pick any $\bar{x} \in E_+$ and any nonnegative $u_{\max}$. Then, using estimate (14),

$$\begin{aligned}
\nabla V(z)\, f^*(z, u) &= \langle \rho(z) - \rho(\bar{x}), f(z) \rangle + \langle \rho(z) - \rho(\bar{x}), C'(u - h(z)) \rangle \\
&\leq -c_1 |D_0(\rho(z) - \rho(\bar{x}))|^2 + \langle \rho(z) - \rho(\bar{x}), C'(u - h(z)) \rangle \\
&\leq \langle C(\rho(z) - \rho(\bar{x})), u - h(z) \rangle \\
&:= \sum_{i=1}^{p} s_i(z, u, \bar{x})
\end{aligned}$$

where $s_i(z, u, \bar{x}) = (\ln h_i(z) - \ln h_i(\bar{x}))(u_i - h_i(z))$.

Now, let $h_{\max}$ be a constant depending only on $\bar{x}$ such that, $\max_{i=1,\ldots,p} |\ln h_i(\bar{x})| \leq h_{\max}$.

For each fixed $z$, define the following finite, disjoint sets of integers

$$I_+ = I_+(z) = \{i : h_i(z) > 1\}$$
$$I_- = I_-(z) = \{i : h_i(z) \leq 1\}.$$

Clearly $I_+ \cup I_- = \{1, \ldots, p\}$, and for each $i \in I_-$,

$$\begin{aligned}
u_i \ln h_i(z) &\leq 0 \\
|u_i \ln h_i(\bar{x})| &\leq u_{\max} |\ln h_i(\bar{x})| \leq u_{\max} h_{\max} \\
|h_i(z) \ln h_i(z)| &\leq \frac{1}{e} \\
|h_i(z) \ln h_i(\bar{x})| &\leq |\ln h_i(\bar{x})| \leq h_{\max},
\end{aligned}$$

so that, for the corresponding $i$th term in the above sum:

$$\begin{aligned} s_i(z, u, \bar{x}) &= u_i \ln h_i(z) - u_i \ln h_i(\bar{x}) - h_i(z) \ln h_i(z) + h_i(z) \ln h_i(\bar{x}) \\ &\leq 0 + u_{\max} h_{\max} + \frac{1}{e} + h_{\max}. \end{aligned}$$

On the other hand, for each $i \in I_+$, $s_i$ can be decomposed into two terms

$$-(\ln h_i(z) - \ln h_i(\bar{x}))(h_i(z) - h_i(\bar{x})) + (\ln h_i(z) - \ln h_i(\bar{x}))(u_i - h_i(\bar{x}))$$

the first of which is always negative. Since $h_i(z) > 1$ there is a Lipschitz constant $c_L = c_L(\bar{x})$ such that $|\ln h_i(z) - \ln h_i(\bar{x})| \leq c_L |h_i(z) - h_i(\bar{x})|$ for all $z$ such that $h_i(z) > 1$. So we have,

$$s_i(z, u, \bar{x}) \leq \begin{cases} 0, & \text{if } |u_i - h_i(\bar{x})| \leq |h_i(z) - h_i(\bar{x})| \\ 2c_L |u_i - h_i(\bar{x})|^2, & \text{if } |u_i - h_i(\bar{x})| > |h_i(z) - h_i(\bar{x})|. \end{cases}$$

In either case, we may just write

$$s_i(z, u, \bar{x}) \leq 2c_L(u_{\max}^2 + h_{\max}^2)$$

whenever $i \in I_+$. Then, with $c_{\bar{x}, u_{\max}} = p \max\{u_{\max} h_{\max} + 1/e + h_{\max}, 2c_L(u_{\max}^2 + h_{\max}^2)\}$,

$$\nabla V(z)\, f^*(z, u) \leq c_{\bar{x}, u_{\max}}, \quad \forall\, z \in \mathbb{R}_{>0}^n, \quad \forall\, u \in [0, u_{\max}]^p,$$

as we wanted to show. ∎

We also have the following $\mathbb{R}_{>0}^n$-completeness result (notice that Corollary 2.11 would follow immediately from Proposition 2.9 in the case of an input $u : [0, +\infty) \to \mathbb{U}_\theta$):

**Corollary 2.11** Under the assumptions of Proposition 2.9, system (15) is $\mathbb{R}_{>0}^n$-complete with input-value set $\mathbb{R}_{\geq 0}^p$.

*Proof.* Suppose that $u : [0, +\infty) \to [0, u_{\max}]^p$, and pick any initial condition $z(0) \in \mathbb{R}_{>0}^n$. We already know, from Proposition 2.9, that system (15) is $\mathbb{R}_{>0}^n$-invariant with input-value set $\mathbb{R}_{\geq 0}^p$. Suppose that the maximal interval of existence would be $[0, t_{\max})$ with $t_{\max} < +\infty$.

From Lemma 2.10 we have that

$$\frac{d}{dt} V(z(t)) = \nabla V(z)\, f^*(z, u) \leq c_{\bar{x}, u_{\max}}, \quad \forall t.$$

So $V(z(t)) \leq V(z(0)) + c_{\bar{x}, u_{\max}} t_{\max}$. Since $V$ is proper (property (ii) of Definition 2.4), we conclude that $z(t)$ belongs to a compact subset of the state space $\mathbb{R}^n$, a contradiction with $t_{\max} < \infty$. ∎

**Proof of Theorem 3:** Pick any initial states $x(0) \in \mathbb{R}_{>0}^n$ and $z(0) \in \mathbb{R}_{>0}^n$ of the original system (1) and the observer, respectively. We let $w(\cdot) = (x(\cdot), z(\cdot))$ be the maximal trajectory of the composite system

$$\begin{aligned} \dot{x} &= f(x) \\ \dot{z} &= f(z) + C'(h(x) - h(z)), \end{aligned}$$

which we also write as $\dot{w} = g(w)$, with initial condition $(x(0), z(0))$. We need to show that $w(t) = (x(t), z(t))$ is defined for all $t > 0$, and $|z(t) - x(t)| \to 0$ as $t \to +\infty$.

Since we know that $x(t)$ is defined for all $t \geq 0$ and converges to some equilibrium $\bar{x}$ as $t \to +\infty$, we must prove that $z(t)$ is also defined for all $t \geq 0$ and converges to this same $\bar{x}$ as $t \to +\infty$.

Fix $\theta$ to be any fixed constant such that $0 < \theta < 1$ and (since $x(t)$ converges) let $T$ be such that
$$t \geq T \quad \Rightarrow \quad |h_i(x(t)) - h_i(\bar{x})| \leq \frac{\theta}{2} h_i(\bar{x})$$
for all $i = 1, \ldots, p$. Let $\mathbb{U}_\theta$ be the set of vectors $u$ such that $|u_i - h_i(\bar{x})| \leq \theta h_i(\bar{x})/2$. Corollary 2.11 shows that the solution $z(t)$ exists for all $t \geq 0$. Now, Proposition 2.9, applied with $u(t) = h(x(t+T))$ and $\bar{u} = h(\bar{x})$ (note that $h(x(t)) \in \mathbb{U}_\theta$ for all $t \geq T$), shows that, for all $t \geq T$,
$$|z(t) - \bar{x}| \leq \beta(|z(T) - \bar{x}|, t) + \varphi(\|h(x) - h(\bar{x})\|).$$
where $\beta(r, t) \to 0$ for each fixed $r$ and $\varphi(r) \to 0$ as $r \to 0$.

Now, given any $\varepsilon > 0$, let $T_1$ be such that
$$\varphi(\|h(x) - h(\bar{x})\|_1) < \frac{\varepsilon}{2}$$
where $\|h(x) - h(\bar{x})\|_1 = \text{ess. sup.} \{|h(x(t)) - h(\bar{x})| : t \geq T_1\}$ (such $T_1$ exists because
$$|h(x(t)) - h(\bar{x})| \to 0$$
as $t \to +\infty$).

Next, choose $T_2 \geq T_1$ such that
$$\beta(|z(T_1) - \bar{x}|, t) < \frac{\varepsilon}{2}, \quad \forall t \geq T_2.$$
Then, rechoosing $T$ (if necessary) to be larger than $T_2$ we have that, for all $t \geq T$, $|z(t) - \bar{x}| \leq \varepsilon$. Therefore, $z(t) \to \bar{x}$ as $t \to +\infty$ as wanted. ∎

**Remark 2.12** The observer (7) can be slightly modified to the form
$$\dot{z} = f(z) + C'W(h(x) - h(z)) \tag{17}$$
where $W$ is any positive definite, diagonal $p \times p$ matrix. This allows more flexibility when choosing the gain matrix, namely, it is possible to assign different weights to each component of $h(x)$. The conclusions of Proposition 2.9 still hold for (17). Since $W = \text{diag}(w_1, w_2, \ldots, w_p)$ with $w_i > 0$ and $W = S'S$ for $S = \text{diag}(\sqrt{w_1}, \sqrt{w_2}, \ldots, \sqrt{w_p})$, $\mathbb{R}^n_{>0}$-invariance is guaranteed by the same argument given in that Proposition. Expression (16) becomes
$$\begin{aligned} f_k^*(\check{z}, x; u) &= \sum_{i=1}^m \sum_{j \in A_0} a_{ij} z_1^{b_{1j}} \ldots z_{k-1}^{b_{(k-1)j}} z_{k+1}^{b_{(k+1)j}} \ldots z_n^{b_{nj}} b_{ki} \\ &+ \sum_{i=1}^m \sum_{j \in A_+} a_{ij} z_1^{b_{1j}} \ldots z_{k-1}^{b_{(k-1)j}} x^{b_{kj}} z_{k+1}^{b_{(k+1)j}} \ldots z_n^{b_{nj}} (b_{ki} - b_{kj}) \\ &+ \sum_{j=1}^p c_{jk} \sqrt{w_k} [u_j - h_j(\check{z}, x)], \end{aligned}$$
which is still nonegative when $x = 0$.

The function $V$ in (12) is also an ISS-Lyapunov function for (17). To check property (iii) we compute (with $\varrho = \rho(z) - \rho(\bar{x})$ and $0 < \theta < 1$):
$$\begin{aligned} \nabla V(z) \, f^*(z, u) &= \langle \varrho, f(z) \rangle + \langle \varrho, C'S'S(u - h(z)) \rangle \\ &= \langle \varrho, f(z) \rangle + \langle SC\varrho, S(u - h(z)) \rangle \\ &= P(z, \bar{x}) + R(z, u, \bar{x}) \end{aligned}$$

where

$$P(z, \bar{x}) = \langle \varrho, f(z) \rangle - (1 - \theta) \langle SC\varrho, S(h(x) - h(z)) \rangle$$
$$R(z, u, \bar{x}) = \langle SC\varrho, S(u - h(\bar{x})) \rangle - \theta \langle SC\varrho, S(h(x) - h(\bar{x})) \rangle.$$

Using the notation $\sigma = SC\varrho$, $\mu = S(h(z) - h(\bar{x}))$ and $v = S(u - h(\bar{x}))$, steps 1 and 2 of the proof of Lemma 2.8 still hold, and we obtain:

$$R(z, u, \bar{x}) \leq c|v|^2$$
$$P(z, \bar{x}) \leq -\alpha_1(|D_0\varrho|) - \alpha_4(|C\varrho|).$$

The function $\alpha_1 \in \mathcal{K}_\infty$ is as in Lemma 2.8. The function $\alpha_4 \in \mathcal{K}_\infty$ and the constant $c$ are similar to $\alpha_2$ and $c_3$, respectively, of that Lemma: $c = 2c_L/\theta$, where $c_L$ is a Lipschitz constant of the logarithmic function on $[A, +\infty)$, when $A = \min\{\sqrt{w_i}h_i(\bar{x})/2 : i = 1, \ldots, p\}$; the function $\alpha_4$ is obtained by repeatedly applying Lemma A.3 to the expression on the left:

$$c_2 \left( \sum_{i \in I_+} |C_i\varrho|^2 + \sum_{i \in I_-} |C_i\varrho| \left(1 - e^{-|C_i\varrho|}\right) \right) \geq \alpha_4(|C\varrho|),$$

since both functions $a(r) = c_2 r^2$ and $b(r) = c_2 r(1 - e^{-r})$ are of class $\mathcal{K}_\infty$. The multiplicative constant is given by $c_2 = \min\{w_1 h_1(\bar{x}), \ldots, w_p h_p(\bar{x})\}$.

With an ISS estimate for the difference $|z(t) - \bar{x}|$ and the completeness result, the same proof above shows that (17) is an observer for system (1).

### 2.3 A Remark on Observation Noise

The ISS estimate obtained for the observer allows us to conclude that the observer is robust with respect to small observation noise. We sketch this next. If the input to the observer is $u(t) = h(x(t)) + n(t)$ instead of $h(x(t))$, then the same conclusions hold regarding global existence of trajectories (at least, provided that $h(x(t)) + n(t)$ is nonnegative). In addition, for large $t$ (so $h(x(t)) \approx h(\bar{x}) = \bar{u}$), the term $\varphi(\|u - \bar{u}\|)$ in the ISS estimate becomes approximately $\varphi(\|n\|)$; thus one obtains an asymptotic estimate on $z(t)$ which is a $\mathcal{K}$ function of the noise level, and is in particular small when $n$ is small in magnitude.

### 2.4 A Separation Property

A direct proof of Theorem 3 can be given using a Lyapunov function for the composite system (plant and observer), which is entirely analogous to the "separation theorem" proof used to establish the convergence of Luenberger observers. Let us sketch the procedure here.

Consider the error $\epsilon = x - z$ and look at the composite system in new coordinates $x$ and $\epsilon$ instead of $x$ and $z$:

$$\Sigma_\epsilon : \begin{array}{l} \dot{x} = f(x) \\ \dot{\epsilon} = f(x) - f(x - \epsilon) - C'(h(x) - h(x - \epsilon)) . \end{array}$$

Now consider the following function:

$$V_2(x, \epsilon) = V(x - \epsilon) = \sum_{i=1}^{n} \bar{x}_i \left[ \frac{x_i - \epsilon_i}{\bar{x}_i} \ln \frac{x_i - \epsilon_i}{\bar{x}_i} + 1 - \frac{x_i - \epsilon_i}{\bar{x}_i} \right].$$

Using arguments similar to those employed earlier, we know that

$$\frac{d}{dt}V_2(x(t), \epsilon(t)) \leq -\alpha(|x - \epsilon - \bar{x}|) + c_3|h(x) - h(\bar{x})|^2$$

for all $x = x(t)$ such that $|h_k(x) - h_k(\bar{x})| < \frac{h_k(\bar{x})}{2}$ for all $k = 1, \ldots, p$. (here we are taking the number $\theta$ from Proposition 2.9 to be equal to $1/2$.) Since $x(t) \to \bar{x}$ as $t \to +\infty$ (Theorem 1), there exists some $T$ such that this inequality holds for every $t > T$.

Since $\epsilon = x - z$ and we have shown that both solutions $x(t)$ and $z(t)$ exist for all $t \geq 0$, we may conclude that a solution $\epsilon(t)$ to the $\epsilon$-equation also exists for all $t \geq 0$ (for a direct proof of this fact, one may use a similar argument to that of Proposition 2.9).

Recall that we already have a Lyapunov function, $V_1(x) = \sum_{i=1}^n \bar{x}_i\, g(x_i/\bar{x}_i)$, for the system $\dot{x} = f(x)$, which satisfies an estimate of the form (see [20]):

$$\frac{d}{dt}V_1(x(t)) \leq -c(\bar{x})\,\delta(x, \bar{x})$$

where

$$\delta(x, \bar{x}) = c_1|x - \bar{x}|^2 + d(|x - \bar{x}|^2)$$

and $d(r) = o(r)$ as $r \to 0$. Pick $\varepsilon < |\bar{x}|/2$ such that

$$|d(|x - \bar{x}|^2)| < \frac{c_1}{2}|x - \bar{x}|^2 \quad \forall\ |x - \bar{x}| < \varepsilon.$$

Now choose the time $T = T(\varepsilon)$ so large that

$$|x(t) - \bar{x}| < \varepsilon \quad \text{and} \quad |h(x(t)) - h(\bar{x})| < \frac{h(\bar{x})}{2}, \text{ for all } t > T.$$

So, for $t > T$, we also have that $x(t)$ evolves in a compact set $B(\bar{x}) = \{a \in \mathbb{R}_{>0}^n : |a - \bar{x}| \leq |\bar{x}|/2\}$, and hence there exists a Lipschitz constant, $c_L$, (take, for example, $c_L = \max |\frac{\partial h_k}{\partial x_j}(x)|$ with the maximum taken over $j = 1, \ldots, n$ and $k = 1, \ldots, p$ and $x \in B(\bar{x})$) such that

$$|h(x) - h(\bar{x})| \leq c_L|x - \bar{x}|, \quad \forall\ x \in B(\bar{x}).$$

Thus for $t > T$, the following estimates hold:

$$\frac{d}{dt}V_1(x(t)) \leq -\frac{c_1 c(\bar{x})}{2}|x - \bar{x}|^2,$$
$$\frac{d}{dt}V_2(x(t), \epsilon(t)) \leq -\alpha(|x - \epsilon - \bar{x}|) + c_3 c_L^2|x - \bar{x}|^2.$$

This means that

$$W(x, \epsilon) = \frac{2(c_3 c_L^2 + 2)}{c_1 c(\bar{x})}V_1(x) + V_2(x, \epsilon)$$

is a Lyapunov function for the system $\Sigma_\epsilon$, with

$$\frac{d}{dt}W(x(t), \epsilon(t)) \leq -\alpha(|x - \epsilon - \bar{x}|) - 2|x - \bar{x}|^2,$$

for all $t > T$.

It is not difficult to verify that this estimate can be rewritten as

$$\frac{d}{dt}W(x(t), \epsilon(t)) \leq -\tilde{\alpha}(|\epsilon|^2 + |x - \bar{x}|^2)$$

for some other $\tilde{\alpha} \in \mathcal{K}_\infty$. (First note that, with $\alpha_2(c) = c^2$ and using Lemma A.3, there exists $\alpha_3 \in \mathcal{K}_\infty$ such that

$$\begin{aligned}
\alpha(|x - \epsilon - \bar{x}|) + \alpha_2(|x - \bar{x}|) &\geq \alpha_3(\sqrt{|x - \epsilon - \bar{x}|^2 + |x - \bar{x}|^2}) \\
&\geq \alpha_3\left(\frac{1}{2}|x - \epsilon - \bar{x}| + \frac{1}{2}|x - \bar{x}|\right) \\
&\geq \alpha_3\left(\frac{1}{2}|\epsilon|\right).
\end{aligned}$$

And then, using Lemma A.3 again, $\alpha_3(|\epsilon|/2) + \alpha_2(|x - \bar{x}|) \geq \alpha_3(\sqrt{|\epsilon|^2 + |x - \bar{x}|^2})$.)

Next, notice that $W(x, \epsilon) \leq \frac{2(c_3 c_L^2 + 2)}{c_1 c(\bar{x})} \nu_2(|x - \bar{x}|) + \nu_2(|x - \epsilon - \bar{x}|)$, where $\nu_2 \in \mathcal{K}_\infty$ is as in (11), so immediately

$$W(x, \epsilon) \leq \kappa \nu_2(2\sqrt{|\epsilon|^2 + |x - \bar{x}|^2}) := \tilde{\nu}_2(|\epsilon|^2 + |x - \bar{x}|^2),$$

where $\kappa = 1 + \frac{2(c_3 c_L^2 + 2)}{c_1 c(\bar{x})}$. Then

$$\frac{d}{dt} W(x(t), \epsilon(t)) \leq -\tilde{\alpha}(\tilde{\nu}_2^{-1}(W(x, \epsilon)))$$

which implies that $W(x(t), \epsilon(t)) \to 0$ as $t \to +\infty$ and so $(x(t), \epsilon(t)) \to (\bar{x}, 0)$, as wanted.

## 3  An Alternative Observer

In this section we construct an alternative observer, which works under the same necessary and sufficient conditions for detectability. As mentioned in the Introduction, this observer exhibits slower convergence but this is traded for more robustness to additive disturbances. Among other interesting features, the observer state never becomes negative, even under arbitrary disturbances. which is a natural requirement on physical grounds (since the state being observed is always nonnegative). This nonnegativity is enforced by a "barrier feedback" reminiscent of optimization techniques.

We now use logarithms of the outputs as input to the observer.

**Theorem 4** Consider system (1) and assume that it is detectable. Then the following system, with state-space $\mathcal{X} = \mathbb{R}_{>0}^n$, is an observer for system (1):

$$\dot{z} = f(z) + C'(\rho(h(x)) - \rho(h(z))). \tag{18}$$

We restricted this to a system with state-space $\mathbb{R}_{>0}^n$ so that the logarithms are well-defined.

Using the convenient notation

$$H(x) := [H_1(x)\ H_2(x) \ldots H_p(x)]' = C\rho(x) \tag{19}$$

where $H_i(x) = c_{i1} \ln x_1 + \cdots + c_{in} \ln x_n$, the observer can be written

$$\dot{z} = f(z) + C'(H(x) - H(z)).$$

**Example 3.1** Consider again the system in Example (2.1). Suppose that the output is given by $h(x) = (x_1 x_2^2, x_1 x_4)'$. Then $H(x) = (\ln x_1 + 2 \ln x_2, \ln x_1 + \ln x_4)$ and we can construct the following observer:

$$\begin{aligned}
\dot{z}_1 &= -kz_1 z_2 + k_3 z_3 + z_4 z_4 &&+ (\ln x_1 + 2\ln x_2 - \ln z_1 - 2\ln z_2) \\
& && + (\ln x_1 + \ln x_4 - \ln z_1 - \ln z_4) \\
\dot{z}_2 &= -kz_1 z_2 + k_3 z_3 + z_4 z_4 &&+ 2(\ln x_1 + 2\ln x_2 - \ln z_1 - 2\ln z_2) \\
\dot{z}_3 &= kz_1 z_2 - (k_3 + \beta_3) z_3 && \\
\dot{z}_4 &= \beta_3 z_3 - k_4 z_4 &&+ (\ln x_1 + \ln x_4 - \ln z_1 - \ln z_4).
\end{aligned}$$

The proof of this theorem involves a technique similar to the one used in the proof of Theorem 3, by studying the stability properties of the more general system

$$\dot{z} = f^*(z, u) := f(z) + C'(u - H(z)),$$

and showing that it satisfies an "input to state stability" condition. The observer for (1) is obtained by letting the input be $u(t) = H(x(t))$.

**Lemma 3.2** Let $f^*(z, u) = f(z) + C'(u - H(\bar{x}))$ and assume $C$ is such that $\mathcal{D} + \operatorname{im} C' = \mathbb{R}^n$. Let $\bar{x}$ denote any point in $E_+$ and $V$ be the function defined in (12). Then there exist class $\mathcal{K}_\infty$ functions $\alpha = \alpha_{\bar{x}}$ and $\gamma = \gamma_{\bar{x}}$ such that

$$\nabla V(z)\, f^*(z, u) \leq -\alpha(|z - \bar{x}|) + \gamma(|u - H(\bar{x})|),$$

for all $z \in \mathbb{R}^n_{>0}$ and all $u \in \mathbb{R}^p$. Moreover, one may pick $\gamma(r) = r^2/2$.

*Proof.* Pick any $\bar{x} \in E_+$. We have (using the notation $\varrho = \rho(z) - \rho(\bar{x})$):

$$\begin{aligned}
\nabla V(z)\, f^*(z, u) &= \langle \rho(z) - \rho(\bar{x}), f(z)\rangle + \langle \rho(z) - \rho(\bar{x}), C'(u - H(z))\rangle \\
&= \langle \varrho, f(z)\rangle + \langle C\varrho, u - H(z)\rangle \\
&= \langle \varrho, f(z)\rangle + \langle C\varrho, u - H(\bar{x})\rangle - \langle C\varrho, C\varrho\rangle
\end{aligned}$$

where the definition $H(z) = C\rho(z)$ was used, as well as the equality $u - H(z) = u - H(\bar{x}) + H(\bar{x}) - H(z)$.

Now, using result (14) and then Cauchy-Schwartz and the inequality $ab \leq \frac{1}{2}(a^2 + b^2)$:

$$\begin{aligned}
\nabla V(z)\, f^*(z, u) &\leq -c(\bar{x})|D_0\, \varrho|^2 - |C\, \varrho|^2 + \langle C\varrho, u - H(\bar{x})\rangle \\
&\leq -c(\bar{x})|D_0\, \varrho|^2 - \frac{1}{2}|C\, \varrho|^2 + \frac{1}{2}|u - H(\bar{x})|^2
\end{aligned}$$

for all $z \in \mathbb{R}^n_{>0}$ and all $u \in \mathbb{R}^p$, where $D_0$ is an $(m-1) \times n$ matrix whose rows (transposed) constitute a basis of $\mathcal{D}$, and $c(\cdot)$ is a positive function. By assumption, we know that the matrix $D = [D_0'\, C']'$ has rank $n$. So, for all $z \in \mathbb{R}^n_{>0}$ and all $u \in \mathbb{R}^p$,

$$\begin{aligned}
\nabla V(z)\, f^*(z, u) &\leq -\min\left\{c(\bar{x}), \frac{1}{2}\right\} |D\, \varrho|^2 + \frac{1}{2}|u - H(\bar{x})|^2 \\
&\leq -\alpha(|z - \bar{x}|) + \gamma(|u - H(\bar{x})|),
\end{aligned}$$

where $\alpha(r) = \min\{c(\bar{x}), 1/2\}\tilde{\alpha}(r)$ and $\tilde{\alpha}$ is as in Lemma A.2 (with $S = \{\bar{x}\}$), and $\gamma(r) = \frac{1}{2}r^2$. ∎

We will consider a system with inputs

$$\dot{z} = f^*(z, u) := f(z) + C'(u - H(z)). \tag{20}$$

Since
$$C'(u - H(z)) = \sum_{k=1}^{n} \sum_{i=1}^{p} c_{ik}(u_i - H_i(z))e_k,$$

we may introduce the set of indices

$$K = \{k : \exists i, \ 1 \le i \le p \text{ with } c_{ik} \ne 0\},$$

which is the complement of all coordinates for which $f_k^*(z, u) \equiv f_k(z)$ (i.e., all $k$ such that the term $\ln z_k$ appears in $f_k^*$ are contained in $K$). When considering the system (20), we view it as evolving on the following state-space:

$$\mathcal{X} = \mathcal{X}_1 \times \cdots \times \mathcal{X}_n$$

where $\mathcal{X}_k = \mathbb{R}$ if $k \notin K$ and $\mathcal{X}_k = (0, +\infty)$ if $k \in K$. Observe that "the term $\ln z_i$ appears anywhere in $f^*$ if and only if it appears in $f_i^*$", so the largest open set where $f^*$ is well defined is indeed $\mathcal{X}$.

**Proposition 3.3** Assume system (1), defined by (2) and (3), to be detectable. Then the system (20) is $\mathbb{R}^n_{>0}$-complete with input-value set $\mathbb{R}^p$.

Furthermore, pick any $\bar{x} \in E_+$ and let $\kappa$ be a positive constant. Define $\mathbb{U}_\kappa$ to be the subset of $\mathbb{R}^p$ consisting of vectors $u$ which satisfy $|u| \le \kappa$. Then, the system (20) is ISS with input-value set $\mathbb{U}_\kappa$ (with respect to $\bar{x}$ and the input $\bar{u} = H(\bar{x})$).

*Proof.* Note that, for each input $u : [0, +\infty) \to \mathbb{R}^p$, $f^*(\cdot, u)$ is of class $C^1$ and locally Lipschitz on $\mathcal{X}$ for each fixed $u$, and $f^*(z, \cdot)$ is measurable and locally integrable on $\mathcal{I}$ for each fixed $z$. Then, given an initial condition $z_0$ and a fixed input $u$ as above, a unique maximal solution of (20) exists on a maximal interval, $J$. Observe that $f^*$ does not extend to $z_k = 0$ and $k \in K$, but while $z_k(t) \in \mathbb{R}^n_{>0}$ no difficulties exist. Then $J$ is not empty, and we will show that $J = [0, +\infty)$.

*Step 1.* To prove $\mathbb{R}^n_{>0}$-invariance (with input-value set $\mathbb{R}^p$), we only need to show that, for each coordinate $k$ with $k \notin K$,

$$z_k(t) > 0, \quad \forall \ t \in J,$$

since all the others satisfy this condition by definition of $\mathcal{X}$.

Assume that one of the coordinates is $\le 0$ at some instant and define

$$t_0 = \inf\{t \in J : z_k(t) = 0 \text{ for some } k \notin K\}.$$

Pick any coordinate $k$ such that $z_k(t_0) = 0$.

Note that for every $k \notin K$, $f_k^*(z, u) \equiv f_k(z)$, so using the notation

$$(\check{z}(t), x) = (z_1(t), \ldots, z_{k-1}(t), x, z_{k+1}(t), \ldots, z_n(t)),$$

we can construct a new function $F : [0, \infty) \times \mathbb{R} \to \mathbb{R}$:

$$F(t, x) = \begin{cases} f_k(\check{z}(t), x; u(t)), & t \in [0, t_0) \\ f_k(\check{z}(t_0), x; u_0), & t \in [t_0, +\infty) \end{cases}$$

(arbitrary $u_0$ in $\mathbb{U}_\kappa$). We have that standard existence and uniqueness results apply to $\dot{x} = F(t, x)$.

*Claim.* $F(t, 0) \ge 0$ for almost all $t \in \mathcal{I}$.

To prove this, write

$$f_k(\check{z}, x; u) = \sum_{i=1}^{m} \sum_{j \in A_0} a_{ij} z_1^{b_{1j}} \ldots z_{k-1}^{b_{(k-1)j}} z_{k+1}^{b_{(k+1)j}} \ldots z_n^{b_{nj}} b_{ki}$$

$$+ \sum_{i=1}^{m} \sum_{j \in A_+} a_{ij} z_1^{b_{1j}} \ldots z_{k-1}^{b_{(k-1)j}} x^{b_{kj}} z_{k+1}^{b_{(k+1)j}} \ldots z_n^{b_{nj}} (b_{ki} - b_{kj})$$

where $A_0 = \{j : b_{kj} = 0\}$ and $A_+ = \{j : b_{kj} > 0\}$.

For $x = 0$ and $t \in \mathcal{I}$:

(a) the second term is zero (since $x = 0$);

(b) the first term is nonnegative since, by definition of $t_0$, $z_i(t) \geq 0$, for all $i \neq k$ and $t \leq t_0$, and either $F(t, 0) = f_k(\check{z}(t), 0; u(t))$ for $t < t_0$, or $F(t, 0) = f_k(\check{z}(t_0), 0; u(t_0))$ for $t \geq t_0$.

This proves the claim.

Now, notice that the initial value problem

$$\dot{x} = F(t, x)$$
$$x(0) = z_k(0)$$

has a solution $x(t)$ in some open interval, $\tilde{J}$, which contains $[0, t_0]$, and moreover, for all $t < t_0$, that solution coincides exactly with the trajectory of the $k$-coordinate: $x(t) \equiv z_k(t)$. Also, by continuity of both $x(t)$ and $z_k(t)$, it must be that $x(t_0) = z_k(t_0)$.

From Lemma A.4, with $\mathcal{X} = \mathbb{R}$ and $J_{x_0} = \tilde{J}$ it follows that $x(t) > 0$ for all $t \in \tilde{J}$, and hence $z_k(t) > x(t)$ for all $t \leq t_0$, which contradicts the fact that $z_k(t_0) = 0$.

To conclude, $z_k(t) > 0$ for all $k \notin K$ and all $t \in J$, and thus the solution $z(t)$ never leaves the set $\mathbb{R}^n_{>0}$ in the maximal interval $J$, proving $\mathbb{R}^n_{>0}$-invariance with input-value set $\mathbb{R}^p$.

*Step 2.* Now assume (to get a contradiction) that the interval $J$ is finite, $J = [0, t_{\max})$, $t_{\max} < +\infty$. From Lemma 3.2, and from boundedness of $u$, we know that the derivative of $V(z(t))$ satisfies

$$\frac{d}{dt} V(z(t)) \leq |u - H(\bar{x})|^2 \leq \kappa^2 + |H(\bar{x})|^2 = c,$$

so

$$V(z(t)) \leq V(0) + \int_0^t c \, dt \leq V(0) + c \, t_{\max}.$$

Since $V^{-1}([0, V(0) + c \, t_{\max}])$ is compact, this inequality implies that $z(t)$ remains bounded for all $t \in J$. (This does not yet establish that $J = [0, +\infty)$, since our state-space is $\mathcal{X}$, not $\mathbb{R}^n$; trajectories could conceivably approach the boundary of the positive orthant for some coordinates in $K$.)

Let $b \in \mathbb{R}$ be a (strictly) positive constant such that $z(t) \in [0, b]^n$ for all $t \in J$. Choose strictly positive constants $\varepsilon_k$ for each $k \in K$ such that:

$$f_k^*(z, u) = f_k(z) + \sum_{i=1}^{p} c_{ik} u_i - \sum_{i=1}^{p} \sum_{j \neq k} c_{ik} c_{ij} \ln z_j - \sum_{i=1}^{p} c_{ik}^2 \ln z_k > 0,$$

whenever each $u_i$ takes values in $[-\|u\|, \|u\|]$, $z$ is a vector in $(0, b]^n$ and $z_k < 2\varepsilon_k$. It is possible to make this choice of $\varepsilon_k$, because we are assuming $c_{ij} \geq 0$ for all $i, j$ and we know that $u$ is bounded and also that $0 < z_j \leq b$ which imply the existence of positive constants $a_i$ such that:

$$|f_k(z)| < a_1, \quad \sum_{i=1}^{p} c_{ik} |u_i| < a_2, \quad \text{and} \quad \sum_{i=1}^{p} \sum_{j \neq k} c_{ik} c_{ij} \ln z_j < a_3.$$

Thus, it is enough to take $\varepsilon_k$ such that

$$\ln(2\varepsilon_k) \leq -\frac{1}{\sum_{i=1}^p c_{ik}^2}(a_1 + a_2 + a_3).$$

(Note that $\sum_{i=1}^p c_{ik}^2 \neq 0$ precisely because $k \in K$.) If the value of $\varepsilon_k$ thus obtained is greater than the initial condition $z_k(0)$, we should rechoose $\varepsilon_k$ to get $\varepsilon_k \leq z_k(0)$.

Now, each coordinate function $z_k(t)$, for $k \in K$ and $t \in J$ satisfies an initial value problem of the form

$$\begin{aligned} \dot{x} &= F(t,x) \\ x(0) &= z_k(0) > 0, \end{aligned}$$

where

$$F(t,x) := f_k^*(z_1(t), \ldots, z_{k-1}(t), x, z_{k+1}(t), \ldots, z_n(t), u(t))$$

is locally Lipschitz on $\mathbb{R}_{>0}$ for each fixed $t$ and measurable on $\mathcal{I}$ for each fixed $x$, and moreover $F(t,x) > 0$ whenever $x < 2\varepsilon_k$. Then Lemma A.5 applies and we conclude that

$$z_k(t) > \varepsilon_k, \quad \forall\, t \in J.$$

Putting these together with the conclusion of step 1, we have that trajectories of (20) stay inside the set

$$\mathcal{C} = [\varepsilon_1, b] \times \cdots \times [\varepsilon_n, b].$$

(where $\varepsilon_k = 0$ if $k \notin K$), during the maximal interval $J$.

Since $\mathcal{C}$ is a compact subset of the domain of definition, $\mathcal{X}$, and since $f^*(z,u)$ (when seen as an explicit function of time through $u = u(t)$), is defined for all $t \in [0, +\infty) \times \mathcal{X}$, it follows by standard ODE arguments (e.g., Proposition C.3.6 in [17]) that the maximal interval is $J = [0, +\infty) \cap \mathcal{I} \equiv [0, +\infty)$, contradicting the existence of a finite $t_{\max}$. Hence, $J$ must be the nonnegative half-line, proving $\mathbb{R}_{>0}^n$-completeness with input-value set $\mathbb{U}_\kappa$.

*Step 3.* Finally, notice that the function (12) is in fact an ISS-Lyapunov function with respect to the point $\bar{x}$ and the input $H(\bar{x})$ for the system (20), because all three properties of Definition 2.4 are satisfied ((i) and (ii) we have previously shown and (iii) follows immediately from Lemma 3.2).

It follows from Lemma 2.7 that (20) is ISS with input-value set $\mathbb{U}_\kappa$ (with respect to the point $\bar{x}$ and the input $H(\bar{x})$). ∎

**Remark 3.4** Proposition 3.3 showed that the given system is ISS with respect to bounded inputs, for each possible bound on inputs. It would be more interesting (although, not needed for our purposes) to know if the same system is ISS with respect to arbitrary inputs. That is, it might be possible to pick the comparison functions in the ISS definition in a manner independent of the bound $\kappa$. It is remarked in [1] that "input semiglobal" ISS (ISS with respect to bounded inputs) indeed implies ISS. However, the proof given in that paper is for the more standard notion of ISS, for systems evolving in $\mathbb{R}^n$, while we have generalized the notion to deal with systems with positive states. Thus the proof does not apply.

## 3.1 Proof of Theorem 4

Pick any initial states $x(0) \in \mathbb{R}^n_{>0}$ and $z(0) \in \mathbb{R}^n_{>0}$ of the original system (1) and the observer, respectively. We let $w(\cdot) = (x(\cdot), z(\cdot))$ be the maximal trajectory of the composite system

$$\begin{aligned} \dot{x} &= f(x) \\ \dot{z} &= f(z) + C'(H(x) - H(z)), \end{aligned}$$

which we also write as $\dot{w} = g(w)$, with initial condition $(x(0), z(0))$. We need to show that $w(t) = (x(t), z(t))$ is defined for all $t > 0$, and $|z(t) - x(t)| \to 0$ as $t \to +\infty$.

The rest of the argument is completely analogous to what was done for the other observer. Since we know that $x(t)$ is defined for all $t \geq 0$ and converges to some equilibrium $\bar{x}$ as $t \to +\infty$, we must prove that $z(t)$ is also defined for all $t \geq 0$ and converges to this same $\bar{x}$ as $t \to +\infty$.

Since $x(t)$ converges to $\bar{x}$, so does $H(x(t))$ converge to $H(\bar{x})$. Let $\kappa$ be a constant such that $|H(x(t))| \leq \kappa$ for all $t$. Let $\mathbb{U}_\kappa$ be the set of vectors $u$ with norm less than or equal to $\kappa$, so that $H(x(t)) \in \mathbb{U}_\kappa$ for every $t \geq 0$. Then Proposition 3.3, applied with $u(t) = H(x(t))$ and $\bar{u} = H(\bar{x})$, shows that the solution $z(t)$ exists for all $t \geq 0$ and that $z(t)$ satisfies

$$|z(t) - \bar{x}| \leq \beta(|z(0) - \bar{x}|, t) + \varphi(\|H(x) - H(\bar{x})\|).$$

where $\beta \in \mathcal{KL}$, and $\varphi \in \mathcal{K}_\infty$.

Now convergence is established in a routine ISS manner: given any $\varepsilon > 0$, pick $T_1$ such that

$$\varphi(\|H(x) - H(\bar{x})\|_1) < \frac{\varepsilon}{2}$$

where $\|H(x) - H(\bar{x})\|_1 = \text{ess.sup}.\{|H(x(t)) - H(\bar{x})| : t \geq T_1\}$. Next, pick $T_2 \geq T_1$ such that

$$\beta(|z(T_1) - \bar{x}|, t) < \frac{\varepsilon}{2}, \quad \forall t \geq T_2.$$

Then, for all $t \geq T_2$, $|z(t) - \bar{x}| \leq \varepsilon$. Therefore, $z(t) \to \bar{x}$ as $t \to +\infty$ as wanted.

## 3.2 Robustness of the Alternative Observer

In Sections 2-3, we have considered the general models with inputs

$$\dot{z} = f(z) + C'(u - h(z)) \tag{21}$$

and

$$\dot{z} = f(z) + C'(u - H(z)) \tag{22}$$

which (besides providing us with observers when $u(t) \equiv h(x(t))$ or $u(t) \equiv H(x(t))$ may be useful in questions regarding the control of the class of chemical reactions described in the Introduction.

We have shown in Proposition 3.3 that, as long as $u$ remains in a set $\mathbb{U}_\kappa$ for all times, trajectories of (22) are well defined for all $t$ and satisfy an input-to-state stability condition.

The system (21) also exhibits an ISS stability property, but only with respect to a more restricted set $\mathbb{U}$, namely, $\|u\|$ has to be uniformly bounded by a fixed constant depending on $\bar{x}$ which, in particular, requires every coordinate of $u$ to be nonnegative.

In fact, if it happens that $u(t)$ becomes less than an (arbitrary) amount $-\varepsilon$ on some interval of time (in the case of a systematic negative disturbance or some negative offset in a control, for instance due to leakage outflows in a chemical reactor), then solutions to (21) may blow up in finite time. This problem can be seen in the following simple example. Take the same network as in Example 3.1 with $n = 3, m = 2$, $\mathcal{D} = \text{span}\{(1, 1, -1)'\}$ and rate constants equal to 1, for simplicity. Let

$$C = \begin{bmatrix} 2 & 0 & 0 \\ 0 & 0 & 2 \end{bmatrix}$$

and consider an input $u(t) = (u_1(t), u_2(t))'$ with $2(u_1(t) + u_2(t)) < -\varepsilon$ for some $\varepsilon > 0$ and all $t$. Then (21) becomes

$$\begin{aligned}
\dot{x}_1 &= -x_1 x_2 + x_3 \quad +2(u_1 - x_1^2) \\
\dot{x}_2 &= -x_1 x_2 + x_3 \\
\dot{x}_3 &= +x_1 x_2 - x_3 \quad +2(u_2 - x_3^2).
\end{aligned}$$

Now look at the variable $x_1 + x_3$:

$$\dot{x}_1 + \dot{x}_3 = 2(u_1 + u_2) - (x_1^2 + x_3^2) \le -\varepsilon - \frac{1}{2}(x_1 + x_3)^2.$$

Viewing this as a differential inequality, we may conclude that $(x_1 + x_3)(t) \le w(t)$ where $w(t)$ is the solution to $\dot{w} = -\varepsilon - \frac{1}{2}w^2$, $w(0) = (x_1 + x_3)(0)$, i.e.,

$$(x_1 + x_3)(t) \le -\sqrt{\varepsilon} \tan\left(\frac{\sqrt{\varepsilon}}{2}t - a_0\right)$$

where $a_0 = \arctan\left(\frac{(x_1+x_3)(0)}{\sqrt{\varepsilon}}\right) > 0$. Therefore, for a suitably chosen initial condition, $w(t)$ blows up to $-\infty$ as $t$ tends to $\frac{\sqrt{\varepsilon}}{2}(a_0 + \frac{\pi}{2})$. Hence, since $(x_1 + x_3)(t) \le w(t)$, $x_1 + x_3$ also has to blow up to $-\infty$ in finite time, and this can happen even if all states started positive.

In contrast, model (22), allows each coordinate of $u$ to have negative values, for an arbitrarily large interval of time, and only requires that $|u(t)| < \kappa$ for all $t$, for some (fixed but arbitratry) positive constant, $\kappa$.

## 4  Generalization to Systems with "Multiple Linkage Classes"

We now look at the more general case of a system with vector field of the form (2) but where the matrix $A = \text{diag}(A_1, \ldots, A_L)$ is block diagonal, and each $A_s$ (of size $m_s$) is itself irreducible and has nonnegative entries (or, at least, there exists a permutation matrix, $P$, such that $PAP^{-1}$ has that diagonal form). The matrix $B$ is also partitioned into $B = [B_1 \cdots B_L]$, where each $B_s$ is of dimension $n \times m_s$ and, since the assumption that $B$ has full rank is still valid, each $B_s$ itself has full rank, $m_s$ ($m_1 + \cdots + m_L = m$). The system (1) can be written as $\dot{x} = f_1(x) + \cdots + f_L(x)$ where each $f_s(x)$ is computed according to formula (2) using $A_s$ and $B_s$.

The number $L$ is called the number of "linkage classes" and denotes the (smallest) number of connected components of the incidence graph $G(A)$. $G(A)$ is the graph whose nodes are the integers $\{1, \ldots, m\}$ and for which there is an edge $j \to i$, iff $a_{ij} > 0$.

To each connected component there corresponds a space

$$\mathcal{D}_s = \text{span } \{b_i - b_j : b_i, b_j \text{ are columns of } B_s\}.$$

The assumptions on the $B_s$ imply that each space $\mathcal{D}_s$ has dimension $m_s - 1$,

$$\mathcal{D} = \mathcal{D}_1 \oplus \cdots \oplus \mathcal{D}_L, \tag{23}$$

(direct sum), and so $\dim \mathcal{D} = (m_1 - 1) + \cdots + (m_L - 1) = m - L$.

**Example 4.1** To illustrate this structure, consider a general enzymatic mechanism with uncompetitive inhibitor, consisting of one enzyme $E$, one substrate $S$, one product $P$ and an uncompetitive inhibitor $I$ ($Q$ and $R$ are intermediate complexes):

$$S + E \xrightarrow{k_1} Q, \quad Q \xrightarrow{k_{-1}} S + E,$$
$$Q \xrightarrow{k_2} P + E, \quad P + E \xrightarrow{k_{-2}} Q,$$
$$Q + I \xrightarrow{k_3} R, \quad R \xrightarrow{k_{-3}} Q + I.$$

There are two linkage classes, $L = 2$:

(i) the first class consisting of the complexes $S + E$, $P + E$ and $Q$;

(ii) the second class consisting of the complexes $Q + I$ and $R$.

For class (i), $S + E$, $P + E$ and $Q$ are the three nodes of $G(A_1)$ and, due to the reversibility of the reactions, it is possible to "connect" any two of these nodes through a path in the graph. The same is true for class (ii).

Notice that the same complex, e.g., $Q$ in the above example, may belong to different connected components (otherwise the problem could be reduced to two completely independent "single linkage" problems).

In [16] it is shown that this system does not admit boundary equilibria in any positive class.

Let $x = (S, P, Q, R, E, I)'$. Then $B = [B_1 \ B_2]$ and $A = \text{diag}(A_1, A_2)$:

$$B_1 = \begin{bmatrix} 1 & 0 & 0 \\ 0 & 1 & 0 \\ 0 & 0 & 1 \\ 0 & 0 & 0 \\ 1 & 1 & 0 \\ 0 & 0 & 0 \end{bmatrix}, \quad B_2 = \begin{bmatrix} 0 & 0 \\ 0 & 0 \\ 1 & 0 \\ 0 & 1 \\ 0 & 0 \\ 1 & 0 \end{bmatrix}$$

and

$$A_1 = \begin{bmatrix} 0 & 0 & k_{-1} \\ 0 & 0 & k_2 \\ k_1 & k_{-2} & 0 \end{bmatrix}, \quad A_2 = \begin{bmatrix} 0 & k_{-3} \\ k_3 & 0 \end{bmatrix}.$$

The space $\mathcal{D}$ is given by $\mathcal{D}_1 + \mathcal{D}_2$:

$$\mathcal{D}_1 = \text{span}\,\{(1, -1, 0, 0, 0, 0)', (1, 0, -1, 0, 1, 0)')\}$$
$$\mathcal{D}_2 = \text{span}\,\{(0, 0, 1, -1, 0, 1)'\}$$

and the function $f$ is given by $f_1 + f_2$:

$$f_1(x) = \begin{bmatrix} k_{-1}Q - k_1 SE \\ -k_{-2}PE + k_2 Q \\ -k_{-1}Q + k_1 SE + k_{-2}PE - k_2 Q \\ 0 \\ k_{-1}Q - k_1 SE - k_{-2}PE + k_2 Q \\ 0 \end{bmatrix}, \quad f_2(x) = \begin{bmatrix} 0 \\ 0 \\ k_{-3}R - k_3 QI \\ -k_{-3}R + k_3 QI \\ 0 \\ k_{-3}R - k_3 QI \end{bmatrix}.$$

For a general "multiple linkage" system we may consider output maps of the same form as before, i.e., monomials in the state variables, as in (3). These monomials may include any of the variables "$x_i$", and they have the same interpretation as before: either representing the concentration of some of the substances (as in the case of $x_1$), or being proportional to some reaction rate (as in the case of $x_1^3 x_4$, etc.).

The necessary and sufficient detectability condition given in Theorem 2 is still valid (Theorem 1 generalizes, as sketched in [20]). The main fact to verify is that (5) still holds for the general case. But, from [20], we know that for an interior point $\bar{x}$, $f(\bar{x}) = 0$ if and only if $f_s(\bar{x}) = 0$ for each $s = 1, \ldots, L$. (That is, $\bar{x}$ is an equilibrium of the entire system if and only if it is an equilibrium of every system $\dot{x} = f_s(x)$; this nontrivial fact follows from the block irreducibility property.) Then, for each $s$, (5) says that, if $\bar{x} \in E_+$, then, for any $\bar{z} \in \mathbb{R}^n_{>0}$,

$$\rho(\bar{x}) - \rho(\bar{z}) \in \mathcal{D}_s^\perp \quad \Longleftrightarrow \quad \bar{z} \in E_+^s,$$

where $E_+^s$ is the set of interior equilibria of $\dot{x} = f_s(x)$ (so $E_+ = E_+^1 \cap \ldots \cap E_+^L$). Equivalently, if $\bar{x} \in E_+$, then, for any $\bar{z} \in \mathbb{R}^n_{>0}$,

$$\rho(\bar{x}) - \rho(\bar{z}) \in \mathcal{D}_1^\perp \cap \ldots \cap \mathcal{D}_L^\perp = \mathcal{D}^\perp \quad \Longleftrightarrow \quad \bar{z} \in E_+.$$

Thus, Equivalence (5) holds for $L > 1$ as well.

So, a general system $\dot{x} = f_1(x) + \cdots + f_L(x)$, $y = h(x)$ is detectable if and only if the matrix $C$ of the (exponents of the) output map satisfies either condition (c) or (d) in Theorem 2. Note that more "linkage classes" mean more information is needed in order for the system to be detectable. For a $n$-dimensional system, the space $\mathcal{D}$ has dimension $m - L$ and detectability implies that the matrix $C$ must have rank $p = n - (m - L)$. As we have seen, a single linkage class requires $p = n - m + 1$, whereas multiple linkage classes require $p = n - m + L$.

In the example above, for detectability of $\dot{x} = f_1(x) + f_2(x)$, $y = h(x)$, $C$ will need to have rank 3. The following output would be a suitable choice: $h(x) = (S^2Q, RI^2, E)'$.

For a detectable "multiple linkage" system, an observer for $\dot{x} = f_1(x) + \cdots + f_L(x)$ is again of the form

$$\dot{z} = f_1(z) + \cdots + f_L(z) + C'(h(x) - h(z)). \tag{24}$$

To prove convergence of the observer, one may use Proposition 2.9, Lemma 2.8 and Corollary 2.11 as before, after checking some points.

The $\mathbb{R}^n_{>0}$-invariance argument is unchanged since the particular forms of $A$ and $B$ still imply that (16) and the corresponding conclusions hold.

To see that Lemma 2.8 holds, we must analyse the term

$$\nabla V(z) \, f(z) = \nabla V(z) \, f_1(z) + \cdots + \nabla V(z) \, f_L(z).$$

For each $s = 1, \ldots, L$, result (14) from [20] holds, so

$$\begin{aligned}
\nabla V(z) \, f_s(z) &= \langle \rho(z) - \rho(\bar{x}), f_s(z) \rangle \\
&\leq -c_s(\bar{x}) \sum_{i,j \dashv B_s} \langle b_i - b_j, \rho(z) - \rho(\bar{x}) \rangle^2 \\
&\leq -c_s(\bar{x}) |D_{0,s}(\rho(z) - \rho(\bar{x}))|^2,
\end{aligned}$$

where $i,j \dashv B_s$ means that only the columns $b_i, b_j$ of $B_s$ are present in the sum. Each matrix $D_{0,s}$ is such that its rows (when transposed) span the space $\mathcal{D}_s$. Therefore, because of (23), the rows (transposed) of the matrix $D_0 = [D'_{0,1} \cdots D'_{0,L}]'$ span the whole of $\mathcal{D}$ and hence

$$\begin{aligned}
\nabla V(z) \, f(z) &\leq -\sum_{s=1}^L c_s(\bar{x}) \, |D_{0,s}(\rho(z) - \rho(\bar{x}))|^2 \\
&\leq -\min c_s(\bar{x}) \, |D_0(\rho(z) - \rho(\bar{x}))|^2.
\end{aligned}$$

The rest of the proof of Lemma 2.8 is unchanged, since the form of the observer is the same as before. Thus the function $V$ satisfies an estimate

$$\nabla V(z)\ f^*(z, u) \leq -\alpha_1(|D_0(\rho(z) - \rho(\bar{x}))|) - \alpha_2(|C(\rho(z) - \rho(\bar{x}))|) + c_3|u - h(\bar{x})|^2$$

where $\alpha_1(r) = \min c_s(\bar{x})\ r^2$ and $\alpha_2 \in \mathcal{K}_\infty$ and the constant $c_3$ are as in the Lemma. Finally, detectability guarantees that $\mathcal{D} + \text{im}\,C' = \mathbb{R}^n$ (Theorem 2), so $D = [D_0'\ C']'$ has rank $n$ which implies $\nabla V(z)\ f^*(z, u) \leq -\alpha(|z - \bar{x}|) + \gamma(|u - h(\bar{x})|)$ for some functions $\alpha, \gamma \in \mathcal{K}_\infty$.

Corollary 2.11 is still valid, as follows from the estimate $\nabla V(z)\ f(z) \leq 0$ verified above.

## 5 Some Simulations

### 5.1 Robustness of the Observers

To numerically test robustness of our observers, we carried out some simulations to explore their responses in two cases: existence of noise in the output measurements and unknown inputs acting in the system we wish to observe.

As a working example we choose $\dot{x} = f(x)$ to be the network proposed by McKeithan in [14] and took the output to be $h(x) = (x_1 x_2^2, x_1 x_4)'$, as in Example 2.1. The constants were taken to be: $k = 0.5, k_3 = 3, k_4 = 2, \beta_3 = 1$.

In Figure 1, the convergence of both observers for this example is shown.

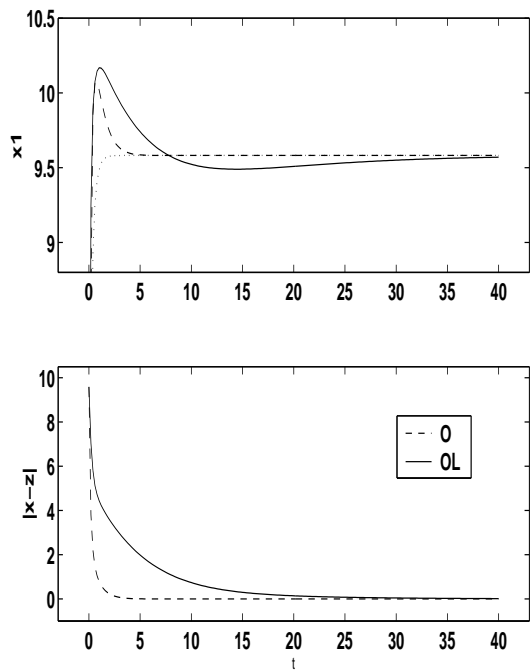

Figure 1: The trajectories of the system and observers without noise and without unknown inputs. The errors, as well as the trajectories of the first coordinate, are shown against time. The dotted line corresponds to the McKeithan network, the dashed line corresponds to our main observer, and the solid line to the alternative observer (where the logarithm of the output is used).

In one simulation, white noise was added to the outputs, so that the equation for the main observer becomes

$$\dot{z} = f(z) + C'(h(x(t)) + n(t) - h(z))$$

and $n(t)$ is an $\mathbb{R}^2$-valued vector white noise. In view of the Section 2.3, solutions of this system exist for all $t \geq 0$, provided that $h(x(t)) + n(t)$ is nonnegative; and the observer should provide estimates close to the true state as long as the magnitude of $n(t)$ is small. Thus we chose $x(0)$ so that $x_i(t) \geq 2$ and $|n(t)| < 2$ for all $t$ and all $i$.

With output noise the alternative observer becomes

$$\dot{z} = f(z) + C'(\rho(h(x(t)) + n(t)) - \rho(h(z)))$$

so the noise appears inside the logarithmic function. But under the conditions above we know that this system has a solution for all $t$ (by Proposition 3.3).

In Figure 2 we can see the trajectories of the first coordinates as a result of this simulation. The alternative observer, although slower in the convergence, has smoothed out the noise effect, whereas the main observer exhibits small magnitude perturbations.

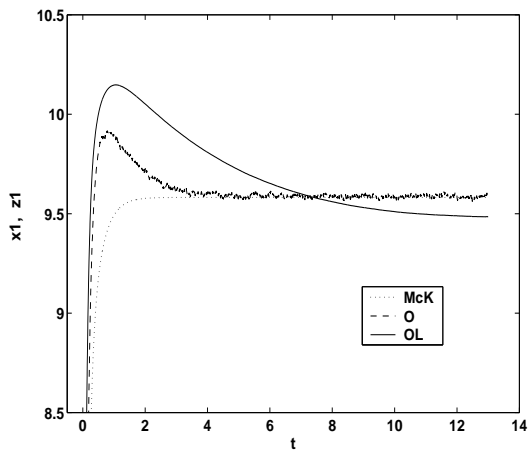

Figure 2: The trajectories of the system and observers in the presence of observation noise.

In another simulation, the model (1) was perturbed by a disturbance consisting of a periodic signal and two "delta" functions. The equation for the model is

$$\dot{x} = f(x) + d(t)$$

where $d(t) = (d_1(t), 0, 0, 0)'$ and $d_1$ is shown in Figure 3. The function $d_1$ was chosen so that (for the same initial condition $x(0)$ as above) $x_i(t) > 0$ for all $i$ and all $t$. (The observers are still the ones for the nominal system, with no disturbance.) Note how our first observer catches-up after the "delta" disturbance, and also tracks (with a small lag) the limit cycle into which the observed system trajectories converge.

In all simulations the initial conditions for system (1) were $x(0) = (3, 2, 3, 20)'$ and for the observers we took $z(0) = (5, 6, 10, 17)'$.

## 5.2 Comparison with Standard Observers

In the chemical reactor literature, observers are typically constructed using an extended Kalman filter (EKF), or, less often, a Luenberger type observer (Lbg) for a linearized system. Neither of these approaches is guaranteed to work for nonlinear systems, but it is often the case that they perform adequately in specific examples, and hence their practical success. In the second part of this section, we compare the performance of our main observer with those of an EKF observer and of a Lbg observer. Our purpose in doing so is to illustrate that, even for very simple examples, these standard techniqes can fail in a major way, while our observer is, as predicted by the theory, convergent to the right estimate.

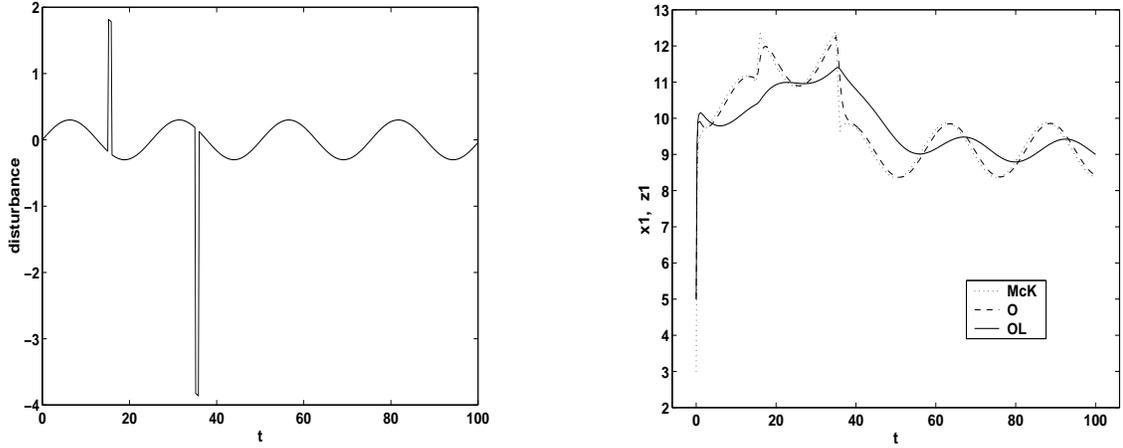

Figure 3: The disturbance shown on the left was added to the system. On the right, its effect on coordinate 1 of the system and observers are plotted against time.

These two standard constructions both have the form

$$\dot{z} = f(z) + L(z)(h(x) - h(z)),$$

where the gain $L(z)$ is to be found in such a way that (at least locally) $|x(t) - z(t)| \to 0$ as $t \to +\infty$. Let us briefly review these constructions.

We consider the linearized dynamics of the error $e = x - z$, around the origin $e = 0$:

$$\dot{e} \doteq [F(z(t)) - L(z(t))H(z(t))]e,$$

where

$$F(z) = Df(e+z)|_{e=0} \quad \text{and} \quad H(z) = Dh(e+z)|_{e=0}$$

are the Jacobians of $f$ and $h$ evaluated at the point $z$.

The gain $L(z)$ for a (continuous) extended Kalman filter is given by

$$L(z(t)) = P(t)H'(z(t))R^{-1},$$

where $P$ is a symmetric positive definite solution to the following Riccati differential equation:

$$\dot{P} = -PH'R^{-1}HP + FP + PF' + Q,$$

and $R$ and $Q$ are two positive definite cost matrices.

A Luenberger type observer is obtained by finding a constant gain $L$ such that the matrix $F(\bar{x}) - LH(\bar{x})$ is Hurwitz. A linearized error equation can also be written as:

$$\dot{e} \doteq [F(x(t)) - LH(x(t))]e$$

(note that the time-dependence of $F$ and $H$ is given in terms of a dependence on the trajectory of the system itself, instead of on the trajectory of the observer). It can be shown that, for initial conditions $x(0)$ and $z(0)$ sufficiently close to $\bar{x}$, this error is asymptotically stable with respect to the origin. (Note that Luenberger observers, at least in their standard formulation, are not a reasonable choice for our example, since their design assumes the knowledge of the equilibrium point around which we are observing. For multi-stable systems such as ours, it makes little sense to assume that this equilibrium is known – in fact, knowing this equilibrium amounts to solving the detectability problem. However, we can still study the behavior of a Luenberger observer, especially since we will show that it does not work even when this additional information is provided.)

Here we take as working example the 2-dimensional system

$$\begin{bmatrix} \dot{x}_1 \\ \dot{x}_2 \end{bmatrix} = \begin{bmatrix} ax_1^2 x_2^5 - x_1^3 x_2^4 \\ -ax_1^2 x_2^5 + x_1^3 x_2^4 \end{bmatrix}, \quad h(x) = x_1^a x_2$$

where $a$ is a positive real number and (in the notation used in the Introduction)

$$A = \begin{bmatrix} 0 & a \\ 1 & 0 \end{bmatrix} \quad \text{and} \quad B = \begin{bmatrix} 3 & 2 \\ 4 & 5 \end{bmatrix}.$$

The invariant classes, $\mathcal{C}_\alpha \subset \mathbb{R}^n_{>0}$, are given by $x_1 + x_2 = \alpha$ and we choose $\alpha = a + 1$, which has the (globally asymptotically stable with respect to the class) equilibrium $\bar{x} = (a, 1)'$.

We would expect both EKF and Lbg to converge provided that $z(0)$ is in a sufficiently small neighborhood of $\bar{x}$, but they may diverge when $z(0)$ is far from $\bar{x}$.

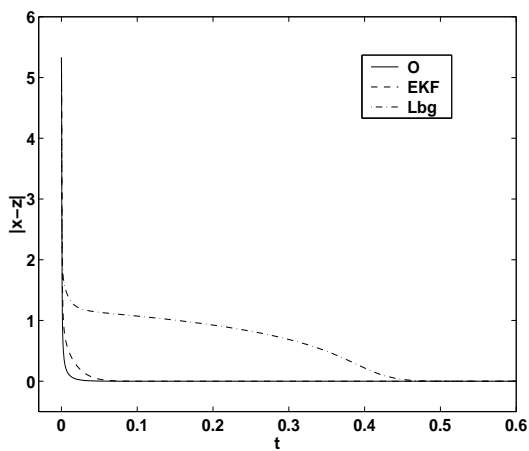

Figure 4: Local convergence of the observers. The initial condition is $z(0) = (4.4, 1.3)'$.

To choose the gain for Lbg we need to pick $L = [l_1 \; l_2]'$ such that the matrix

$$F - LH \mid_{\bar{x}} = \begin{bmatrix} -a^2 - l_1 a^a & a^3 - l_1 a^a \\ a^2 - l_2 a^a & -a^3 - l_2 a^a \end{bmatrix}$$

is Hurwitz (that is, all its eigenvalues have negative real parts). For the simulations shown here let $a = 4$ and $L = [-0.7 \; 1.0]'$.

In Figure 4, it is seen that all observers converge locally, although the errors converge to zero in different time scales.

For other initial conditions, away from $\bar{x}$, both EKF and Lbg seem to diverge, and their performances are clearly inferior to that of our observer. In Figures 5 and 6 the behavior of the three observers is shown in several regions of the positive orthant.

In all simulations the initial conditions for the system were $x(0) = (0.3, 4.7)'$ and the gain matrix for the Luenberger observer was $L = [-0.7 \; 1.0]'$. For the Riccati differential equation the initial condition was $P(0) = 0.1I$, and the cost matrices $Q, R$ were constant and equal to $I$, $I$ being the $2 \times 2$ identity matrix.

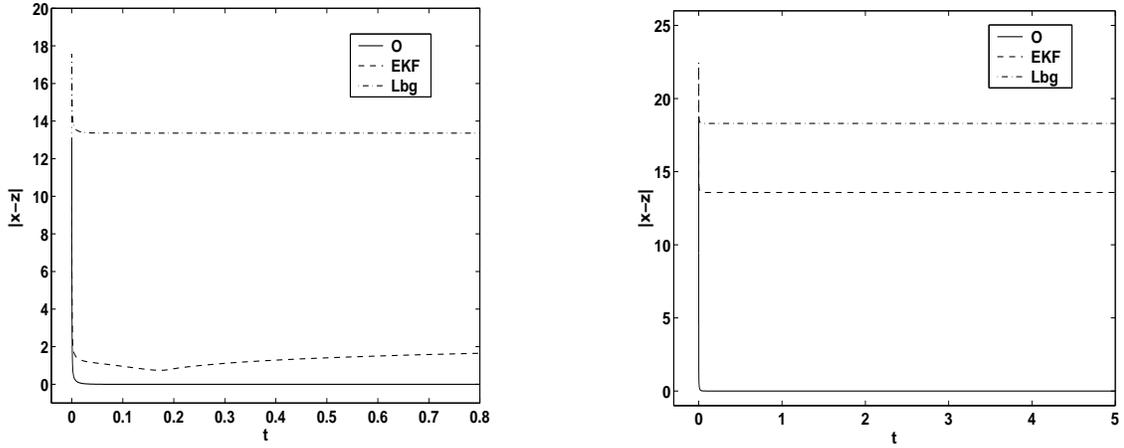

Figure 5: Both extended Kalman filter and Luenberger observer diverge for initial conditions on the higher part of the positive orthant of the plane. On the left plot $z(0) = (3, 16)'$ and on the right plot $z(0) = (19, 6)'$.

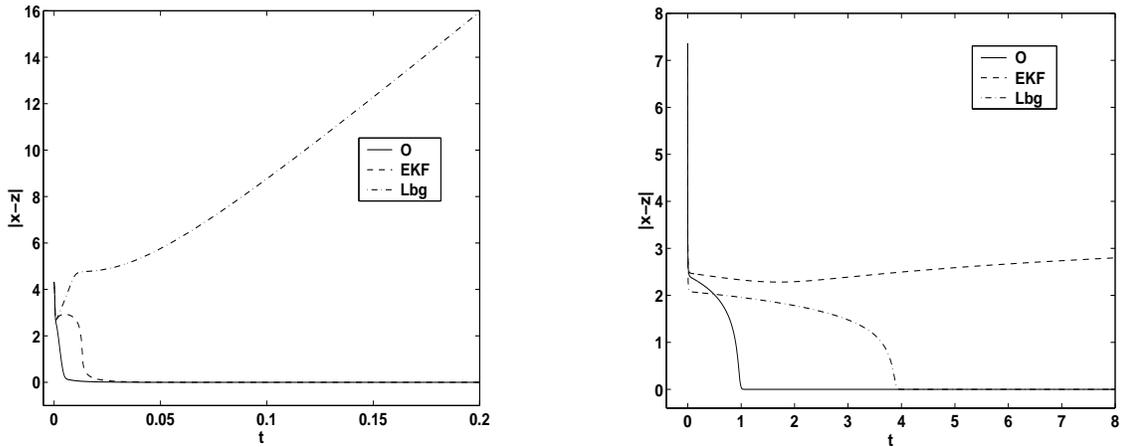

Figure 6: There are regions (close to the x-axis) where one of observer converges but the other diverges. On the left plot $z(0) = (0.8, 0.4)'$ and Lbg diverges. On the right plot $z(0) = (6, 0.1)'$ and EKF diverges.

# A Appendix

## A.1 Upper and Lower Bounds for an ISS-Lyapunov Function

For a general ISS-Lyapunov function defined according to Definition 2.4, one can show the existence of class $\mathcal{K}_\infty$ bounds, as follows. Let

$$\nu_1(r) = \inf\{V(z) : |z - \bar{x}| \geq r \text{ and } z \in \mathbb{R}^n_{\geq 0}\}$$

and

$$\nu_2(r) = \max\{V(z) : |z - \bar{x}| \leq r \text{ and } z \in \mathbb{R}^n_{\geq 0}\}.$$

The function $\nu_2$ is continuous because $V$ is, nondecreasing and satisfies $\nu_2(0) = 0$. Without loss of generality, we can say that it is strictly increasing (by taking $\nu_2(r) + r$).

The function $\nu_1$ is finite, because for any $r > 0$, we may pick any $x_r$ with $|x_r - \bar{x}| = r$ and consider $L := V(x_r)$. Then

$$\nu_1(r) = \min\{V(z) : |z - \bar{x}| \geq r \text{ and } z \in \mathbb{R}^n_{\geq 0} \text{ and } V(z) \leq L\}$$

(the minimum over this set exists, because it is a compact). Thus $\nu_1$ is also continuous and nondecreasing and $\nu_1(0) = 0$. Without loss of generality, $\nu_1$ can be assumed to be strictly increasing.

In the particular case of the function $V$ defined in (12) the functions $\nu_1, \nu_2$ may be taken as follows. Consider $v : \mathbb{R}_{\geq 0} \to \mathbb{R}_{\geq 0}$ given by:

$$v(a) = \begin{cases} (1-a)\ln(1-a) + a, & 0 \leq a \leq 1 \\ (1+a)\ln(1+a) - a + 2(1 - \ln 2), & 1 \leq a < +\infty. \end{cases}$$

Then $v$ is a $\mathcal{K}_\infty$ function (note that $2(1 - \ln 2)$ is a positive quantity that guarantees continuity of $v$ at $a = 1$).

From (12), recall that $g(r) = r \ln r + 1 - r$ for $r \geq 0$. Equivalently,

$$g(1 + a) = (1 + a)\ln(1 + a) - a, \text{ for } a \geq -1.$$

Then

$$g(1+a) \leq \begin{cases} (1-|a|)\ln(1-|a|) + |a|, & |a| \leq 1 \\ (1+a)\ln(1+a) - a + 2(1 - \ln 2), & 1 \leq a < +\infty. \end{cases}$$

So $g(z_i/\bar{x}_i) \leq v(|z_i - \bar{x}_i|/\bar{x}_i)$ and

$$V(z) \leq \sum_{i=1}^n \bar{x}_i\, v\left(\frac{|z_i - \bar{x}_i|}{\bar{x}_i}\right).$$

By using the fact that

$$v\left(\frac{|z_i - \bar{x}_i|}{\bar{x}_i}\right) \leq v\left(\frac{|z_i - \bar{x}_i|}{\min \bar{x}_i}\right),$$

and also keeping in mind that $v(a_i) \leq v(\sqrt{a_1^2 + \cdots + a_n^2})$, it follows that

$$V(z) \leq n \max \bar{x}_i\, v(|z - \bar{x}|/\min \bar{x}_i) = \nu_2(|z - \bar{x}|).$$

For the lower bound, consider

$$w(a) = (1+a)\ln(1+a) - a, \quad 0 \leq a \leq +\infty$$

and notice that

$$g(1+a) \geq \begin{cases} (1+|a|)\ln(1+|a|) - |a|, & |a| \leq 1 \\ (1+a)\ln(1+a) - a, & 1 \leq a < +\infty. \end{cases}$$

So, $g(z_i/\bar{x}_i) \geq w(|z_i - \bar{x}_i|/\bar{x}_i)$ and

$$V(z) \geq \sum_{i=1}^n \bar{x}_i w\left(\frac{|z_i - \bar{x}_i|}{\bar{x}_i}\right).$$

Using Lemma A.3 repeatedly, we may conclude that $V$ is lower bounded by a $\mathcal{K}_\infty$ function.

## A.2 Auxiliary Calculations

The next lemmas state some simple inequalities which are used throughout the text.

**Lemma A.1** Let $v, w \in ]-1, +\infty[$. Then

$$(\ln(1+v))^2 + (\ln(1+w))^2 \geq \frac{1}{2}\left(\ln(1+\sqrt{v^2+w^2})\right)^2.$$

*Proof.* First, note that $(\ln(1+v))^2 \geq (\ln(1-v))^2$ for all $-1 < v < 0$ by putting,

$$g(v) = (\ln(1+v))^2 - (\ln(1-v))^2 = (\ln(1+v) + \ln(1-v))(\ln(1+v) - \ln(1-v)).$$

The second factor is clearly negative and the first factor is strictly increasing on the interval $]-1, 0[$ (by computing the derivative) and is 0 when $v = 0$, therefore it must be also negative on this interval, implying that $g(v) \geq 0$ as desired.

Then it holds that

$$(\ln(1+v))^2 \geq (\ln(1+|v|))^2, \quad \forall v \in ]-1, +\infty[,$$

so it is enough to prove the lemma for all $v, w \in [0, \infty[$. Using the fact that $x^2$ is a convex function, it follows that

$$\frac{1}{2}(\ln(1+v))^2 + \ln(1+w))^2) \geq \frac{1}{4}(\ln(1+v) + \ln(1+w))^2 = \frac{1}{4}(\ln(1+v+w+vw))^2. \quad (25)$$

Moreover, since both $v$ and $w$ are positive: $w + v \geq \sqrt{v^2 + w^2}$. The logarithm is an increasing function, so

$$\frac{1}{4}(\ln(1+v+w+vw))^2 \geq \frac{1}{4}(\ln(1+\sqrt{v^2+w^2}))^2. \quad (26)$$

Putting (25) and (26) together gives desired result. ∎

**Lemma A.2** Let $m \geq n$ and let $D$ be an $m \times n$ matrix of rank $n$, and consider a fixed compact subset $S$ of $\mathbb{R}^n_{>0}$. Then, there exists a class $\mathcal{K}_\infty$ function $\alpha$ such that $|D(\rho(x) - \rho(a))|^2 \geq \alpha(|x-a|)$ for all $a \in S$ and all $x \in \mathbb{R}^n_{>0}$.

*Proof.* Since $D$ has maximal rank, it has a left pseudo-inverse, for example $D^\# = (D'D)^{-1}D'$. Then, for any vector $\varrho \in \mathbb{R}^n$,

$$|\varrho| = |D^\# D\varrho| \leq \|D^\#\| |D\varrho| \quad \Longrightarrow \quad |D\varrho| \geq \frac{1}{\|D^\#\|}|\varrho|.$$

Also, if $\varrho = \rho(x) - \rho(a)$,

$$\varrho_i = \ln\frac{x_i}{a_i} = \ln\left(1 + \frac{x_i - a_i}{a_i}\right) := \ln(1+v_i).$$

Then, applying Lemma A.1 inductively on $n$ we have

$$|\varrho|^2 \geq \frac{1}{2^{n-1}}\left(\ln\left(1+\sqrt{v_1^2 + \cdots + v_n^2}\right)\right)^2 \geq \frac{1}{2^{n-1}}\left(\ln\left(1+\frac{|x-a|}{\max a_i}\right)\right)^2.$$

Now, clearly

$$\alpha(w) = \frac{1}{2^{n-1}\|D^{\#}\|^2}\left(\ln\left(1+\frac{w}{s}\right)\right)^2$$

with $s$ an upper bound on the magnitudes of the coordinates of points of $S$, is of class $\mathcal{K}_\infty$, which proves the lemma.

∎

The following is easy to prove using $\alpha_3(c) = \min\{\alpha_1(c/2), \alpha_2(c/2)\}$:

**Lemma A.3** Assume that $\alpha_1, \alpha_2$ are two $\mathcal{K}_\infty$ functions. Let $a$ and $b$ be positive numbers. Then there exists another $\mathcal{K}_\infty$ function, $\alpha_3$, such that

$$\alpha_1(a) + \alpha_2(b) \geq \alpha_3(\sqrt{a^2+b^2}).$$

for all $a, b > 0$.

## A.3  Some Simple Facts Concerning Invariance

Consider the scalar initial value problem

$$\dot{x} = F(t,x) \tag{27}$$
$$x(0) = x_0$$

where the function $F$ is assumed to have domain $\mathcal{I} \times \mathcal{X}$, where $\mathcal{X}$ is an open subset of $\mathbb{R}$ and $\mathcal{I} = [0, +\infty)$. Let $F$ be locally Lipschitz in $x$ and measurable in $t$, more precisely,

(i) For each $a \in \mathcal{X}$ there exists a real number $r_a$ and a locally integrable function $\alpha : \mathbb{R} \to [0, +\infty)$ such that the ball of radius $r_a$ centered at $a$, $B_{r_a}(a) \subset \mathcal{X}$ and

$$\|F(t,x) - F(t,y)\| \leq \alpha(t)\|x - y\|$$

for each $t \in \mathbb{R}$ and $x, y \in B_{r_a}(a)$.

(ii) For each fixed $a \in \mathcal{X}$, the function $g : \mathcal{I} \to \mathcal{X}$ given by $g(t) := F(t,a)$ is measurable.

For each $x_0 \in \mathcal{X}$ let $J = J_{x_0}$ be the maximal interval of existence of solutions of (27) in forward time. This is an interval of the form $[0, t_{\max})$ with $0 < t_{\max} \leq +\infty$.

Using a standard comparison principle (as done for instance in [20]), we have:

**Lemma A.4** Consider system (27) with domain $\mathcal{X} = \mathbb{R}$ and assume further that,

$$x = 0 \Rightarrow F(t,0) \geq 0 \quad \forall t \in \mathcal{I}.$$

Assume also that the initial condition is positive: $x_0 > 0$. Then $x(t) > 0 \quad \forall t \in J$, i.e., the solution of (27) remains positive for all times in $J$.

**Lemma A.5** Consider again the scalar initial value problem (27), but now with domain $\mathcal{X} = (0, +\infty)$. Assume that the initial condition is positive: $x_0 > 0$, and that the following property holds, for some $0 < \varepsilon < x_0$:

$$\forall t \in \mathbb{R}, \quad F(t,x) > 0 \text{ whenever } x < 2\varepsilon. \tag{28}$$

Then for all $t \in J$, the solutions of (27) satisfy $x(t) > \varepsilon$.

*Proof.* Suppose that the result is false, and let

$$s = \inf_{t \in J}\{t : x(t) \leq \varepsilon\}.$$

We may assume that $s < \infty$ and $s \in J$, because otherwise the desired result is trivial. In other words, let us suppose that the solution $x(t)$ of (27) hits the value $\varepsilon$ at some finite instant $s$ in $J$.

Continuity of solutions gives

$$\exists \delta = \delta(\varepsilon) \quad \text{such that} \quad |x(t) - \varepsilon| < \varepsilon/2 \quad \forall t \in (s - \delta, s + \delta).$$

By property (28), $F(t, x(t)) > 0 \ \forall t \in (s-\delta, s+\delta)$ implying that $x(t)$ is strictly increasing on this $\delta$-interval and so $x(s) > x(s - \delta)$. But, then $x(s - \delta) < \varepsilon$, which contradicts the definition of $s$. ∎